# DEFORMATIONS AND STABILITY IN COMPLEX HYPERBOLIC GEOMETRY

Boris Apanasov†

ABSTRACT. This paper concerns with deformations of noncompact complex hyperbolic manifolds (with locally Bergman metric), varieties of discrete representations of their fundamental groups into $PU(n,1)$ and the problem of (quasiconformal) stability of deformations of such groups and manifolds in the sense of L.Bers and D.Sullivan.

Despite Goldman-Millson-Yue rigidity results for such complex manifolds of infinite volume, we present different classes of such manifolds that allow non-trivial (quasi-Fuchsian) deformations and point out that such flexible manifolds have a common feature being Stein spaces. While deformations of complex surfaces from our first class are induced by quasiconformal homeomorphisms, non-rigid complex surfaces (homotopy equivalent to their complex analytic submanifolds) from another class are quasiconformally unstable, but nevertheless their deformations are induced by homeomorphisms.

1991 *Mathematics Subject Classification.* 57, 55, 53, 51, 32, 22, 20.

†Supported in part by the NSF; research at the University of Tokyo was supported in part by JSPS.

Typeset by $\mathcal{A}_{\mathcal{M}}\mathcal{S}$-TEX

# Deformations and Stability in Complex Hyperbolic Geometry

Boris Apanasov

## 1. Introduction and formulation of main results

This paper presents a recent progress in the theory of deformations of noncompact complex hyperbolic manifolds $M$ (of infinite volume, with variable sectional curvature) and spherical Cauchy-Riemannian manifolds at their infinity $M_\infty$, varieties of discrete faithful representations of the fundamental groups $\pi_1 M$ into $PU(n,1)$, and the problem of (quasiconformal) stability of deformations of such groups and manifolds whose geometry makes them surprisingly different from those in the real hyperbolic geometry with constant negative sectional curvature.

Geometry of the complex hyperbolic space $\mathbb{H}^n_\mathbb{C}$ is the geometry of the unit ball $\mathbb{B}^n_\mathbb{C}$ in $\mathbb{C}^n$ with the Kähler structure given by the Bergman metric whose automorphisms are biholomorphic automorphisms of the ball, i.e., elements of $PU(n,1)$. We notice that complex hyperbolic manifolds (modeled on $\mathbb{H}^n_\mathbb{C}$) with non-elementary fundamental groups are complex hyperbolic in the sense of S.Kobayashi [Kob]; we refer the reader to [AX1, CG, G4] for general information on such manifolds, in particular for several equivalent descriptions of the basic class of geometrically finite complex hyperbolic manifolds and for a discussion on surprising differences between such manifolds and real hyperbolic manifolds with constant negative sectional curvature. Here we study deformations of complex hyperbolic manifolds and their fundamental groups by using the spherical Cauchy-Riemannian geometry at infinity. This CR-geometry is modeled on the one point compactification of the (nilpotent) Heisenberg group $\mathcal{H}_n$, which appears as the sphere at infinity of the complex hyperbolic space $\mathbb{H}^n_\mathbb{C}$. Since any complex hyperbolic manifold can be represented as the quotient $M = \mathbb{H}^n_\mathbb{C}/G$ by a discrete torsion free isometry action of the fundamental group of $M$, $\pi_1(M) \cong G \subset PU(n,1)$, its boundary at infinity $\partial_\infty M$ is naturally identified as the quotient $\Omega(G)/G$ of the discontinuity set of $G$ at infinity. Here the discontinuity set $\Omega(G)$ is the maximal subset of $\partial\mathbb{H}^n_\mathbb{C}$ where $G$ acts discretely; its complement $\Lambda(G) = \partial\mathbb{H}^n_\mathbb{C}\setminus\Omega(G)$ is the limit set of $G$, $\Lambda(G) = \overline{G(x)} \cap \partial\mathbb{H}^n_\mathbb{C}$ for any $x \in \mathbb{H}^n_\mathbb{C}$.

One can reduce the study of deformations of complex hyperbolic manifold $M$, or equivalently the Teichmüller space $\mathcal{T}(M)$ of isotopy classes of complex hyperbolic structures on $M$, to studying the variety $\mathcal{T}(G)$ of conjugacy classes of discrete faithful representations $\rho : G \to PU(n,1)$ (involving the space $\mathcal{D}(M)$ of the developing maps, see [G2, FG]). Here $\mathcal{T}(G) = \mathcal{R}_0(G)/PU(n,1)$, and the variety $\mathcal{R}_0(G) \subset \mathrm{Hom}(G, PU(n,1))$ consists of discrete faithful representations $\rho$ of the group $G$ whose co-volume $\mathrm{Vol}(\mathbb{H}^n_\mathbb{C}/G)$ may be infinite.

Due to the Mostow rigidity theorem [Mo1], hyperbolic structures of finite volume and (real) dimension at least three are uniquely determined by their topology, and one has no continuous deformations of them. Despite that, real hyperbolic manifolds $N$ can be deformed as conformal manifolds, or equivalently as higher-dimensional hyperbolic manifolds $M = N \times (0,1)$ of infinite volume. First such deformations were given by the author [A2] and, after Thurston's "Mickey Mouse" example [T], they were called bendings of $N$ along its totally geodesic hypersurfaces, see also [A1, A3-A5, JM, Ko]. Furthermore such a flexibility of the real hyperbolic geometry is emphasized by the fact that all those deformations can be induced by continuous families of $G$-equivariant quasiconformal self-homeomorphisms $f_t : \overline{\mathbb{H}^n_\mathbb{R}} \to \overline{\mathbb{H}^n_\mathbb{R}}$ of the closure of the real hyperbolic space $\mathbb{H}_{\mathbb{R}}^n$. In particular, these $G$-equivariant quasiconformal homeomorphisms deform continuously the limit set $\Lambda(G) \subset \partial\mathbb{H}_{\mathbb{R}}^n$ (of a "Fuchsian group" $G \subset \operatorname{Isom} \mathbb{H}_{\mathbb{R}}^{n-1} \subset \operatorname{Isom} \mathbb{H}_{\mathbb{R}}^n$) from a round sphere $\partial\mathbb{H}_{\mathbb{R}}^{n-1} = S^{n-2} \subset S^{n-1} = \partial\mathbb{H}_{\mathbb{R}}^n$ into nondifferentiably embedded (nonrectifiable) topological $(n-2)$-spheres in $\partial\mathbb{H}_{\mathbb{R}}^n$ which are the limit sets $\Lambda(G_t)$ of "quasi-Fuchsian groups" $G_t = f_t G f_t^{-1} \subset \operatorname{Isom} \mathbb{H}_{\mathbb{R}}^n$, and obviously the restrictions $f_t|_{\Lambda(G)} : \Lambda(G) \to \Lambda(G_t)$ are quasisymmetric maps.

Such a geometric realization of isomorphisms of discrete groups became the start point in our study of deformations of discrete groups of isometries of negatively curved spaces $X$, see [A7]:

**Problem.** *Given an isomorphism $\varphi : G \to H$ of geometrically finite discrete groups $G, H \subset \operatorname{Isom} X$, find subsets $X_G, X_H \subset \overline{X}$ invariant for the action of groups $G$ and $H$, respectively, and an equivariant homeomorphism:*

$$f_\varphi : X_G \to X_H \quad \varphi(g) \circ f_\varphi = f_\varphi \circ g \quad \text{for all } g \in G\,,$$

*which induces the isomorphism $\varphi$. Determine metric properties of $f_\varphi$, in particular whether it is either quasisymmetric or quasiconformal with respect to the given negatively curved metric $d$ in $X$ (or the induced sub-Riemannian Carnot-Carathéodory structure at infinity $\partial X$).*

If the groups $G, H \subset \operatorname{Isom} X$ are neither lattices nor trivial and have parabolic elements, the only known geometric realization of their isomorphisms in dimension $\dim X \geq 3$ is due to P.Tukia's [Tu] isomorphism theorem for real hyperbolic spaces $X = \mathbb{H}_{\mathbb{R}}^n$. However, the Tukia's construction (based on geometry of convex hulls of the limit sets $\Lambda(G)$ and $\Lambda(H)$) cannot be used in the case of variable negative curvature due to lack of control over convex hulls (where the convex hull of three points may be 4-dimensional), especially nearby parabolic fixed points. However, as a first step in solving the above geometrization Problem, we have the following isomorphism theorem [A8, A10] in the complex hyperbolic space:

**Theorem 3.2.** *Let $\varphi : G \to H$ be a type preserving isomorphism of two non-elementary geometrically finite discrete subgroups $G, H \subset \operatorname{Isom} \mathbb{H}_{\mathbb{C}}^n$. Then there exists a unique equivariant homeomorphism $f_\varphi : \Lambda(G) \to \Lambda(H)$ of their limit sets that induces the isomorphism $\phi$.*

However, in contrast to the real hyperbolic case where such geometric realizations of type preserving isomorphisms of geometrically finite groups are always quasisymmetric maps [Tu], it is doubtful that the (unique) equivariant homeomorphism $f_\varphi : \Lambda(G) \to \Lambda(H)$ constructed in Theorem 3.2 is always CR-quasisymmetric (with respect to the CR-structure on the Heisenberg group $\mathcal{H}_n = \partial\mathbb{H}_{\mathbb{C}}^n\backslash\{\infty\}$). Namely, a possible obstruction to quasisymmetricity directly appear from the following corollary of our construction in Section 4:

**Corollary 4.2.** *Let $M = \mathbb{H}_{\mathbb{C}}^2/G$ be a complex hyperbolic surface with the holonomy group $G \subset PU(1,1) \subset PU(2,1)$ that represents the total space of a non-trivial disk bundle over a Riemann surface of genus $p \geq 0$ with at least four punctures (hyperbolic 2-orbifold with at least four punctures). Then the Teichmüller space $\mathcal{T}(M)$ contains a smooth simple curve $\alpha : [0, \pi/2) \hookrightarrow \mathcal{T}(M)$ with the following properties:*

  (1) *the curve $\alpha$ passes through the surface $M = \alpha(0)$;*

(2) *each complex hyperbolic surface $M_t = \alpha(t) = \mathbb{H}_\mathbb{C}^2/G_t$, $t \in [0, \pi/2)$, with the holonomy group $G_t \subset PU(2,1)$ is homeomorphic to the surface $M$;*
(3) *for any parameter $t$, $0 < t < \pi/2$, the complex hyperbolic surface $M_t$ is not quasiconformally equivalent to the surface $M$.*

Besides the claims in this Corollary, it follows also from the construction of the above complex hyperbolic surfaces $M$ and $M_t$ that their boundaries, the spherical CR-manifolds $N = \partial M = \Omega(G)/G$ and $N_t = \partial M_t = \Omega(G_t)/G_t$ have similar properties. Namely these 3-dimensional CR-manifolds $\{N_t\}$, $0 \leq t < \pi/2$, represent points of a smooth simple curve $\alpha_\infty : [0, \pi/2) \hookrightarrow \mathcal{T}(N)$ in the Teichmüller space $\mathcal{T}(N)$ of the manifold $N = N_0 = \partial M$, are mutually homeomorphic total spaces of non-trivial circle bundles over a Riemann surface of genus $p \geq 0$ with at least four punctures, however none of $\{N_t\}$ with $t > 0$ is quasiconformally equivalent to $N = N_0$. We note that, for the simplest case of manifolds with cyclic fundamental groups, a similar (though based on different ideas) effect of homeomorphic but not quasiconformally equivalent spherical CR-manifolds has been also recently presented by R. Miner [Mi].

It is quite natural that the result in Corollary 4.2 is related to the classical problem of quasiconformal stability of deformations from the theory of Kleinian groups, in particular to well known results by L. Bers [Be1, Be2] and D. Sullivan [Su1]. Following to L. Bers [Be2], a finitely generated Kleinian group $G \subset \mathrm{PSL}(2,\mathbb{C})$ is said to be quasiconformally stable if every homomorphism $\chi : G \to \mathrm{PSL}(2,\mathbb{C})$ preserving the square traces of parabolic and elliptic elements (hence type-preserving) and sufficiently close to the identity is induced by an equivarint quasiconformal mapping $w : \overline{\mathbb{C}} \to \overline{\mathbb{C}}$, $\chi(g) = wgw^{-1}$ for all $g \in G$. Due to a Bers's criterium (involving the quadratic differentials for $G$, see [Be1]), it follows that Fuchsian groups, Schottky groups, groups of Schottky type and ceratin non-degenerate $B$-groups are all quasiconformally stable [Be2]. Changing the condition on homomorphisms $\chi$ in terms of the trace of elements $g \in G$ to the condition that $\chi$ preserves the type of elements of a discrete groups $G$, one has a natural generalization of quasiconformal stability for discrete groups $G \subset PU(n,1)$. In that sense, B. Aebisher and R. Miner [AM] recently proved that (classical) Schottky subgroups $G \subset PU(n,1)$ are quasiconformally stable. Nevertheless, as our construction in Section 4 shows, Fuchsian groups $G \subset PU(2,1)$ are quasiconformally unstable:

**Theorem 4.1.** *There are co-finite Fuchsian groups $G \subset PU(1,1) \subset PU(2,1)$ with signatures $(g, r; m_1, \ldots, m_r)$, where genus $g \geq 0$ and there are at least four cusps (with branching orders $m_i = \infty$), such that:*

(1) *the Teichmüller space $\mathcal{T}(G)$ of such a group $G$ contains a smooth simple curve $\alpha$, $\alpha : [0, \pi/2) \hookrightarrow \mathcal{T}(G)$, that passes through the Fuchsian group $G = \alpha(0)$, and whose points $\alpha(t) = G_t \subset PU(2,1)$, $0 < t < \pi/2$, are all non-trivial quasi-Fuchsian groups;*
(2) *each isomorphism $\chi : G \to G_t$, $0 < t < \pi/2$, is induced by a $G$-equivariant homeomorphism $f_t : \overline{\mathbb{H}_\mathbb{C}^2} \to \overline{\mathbb{H}_\mathbb{C}^2}$ of the closure $\overline{\mathbb{H}_\mathbb{C}^2} = \mathbb{H}_\mathbb{C}^2 \cup \partial \mathbb{H}_\mathbb{C}^2$ of the complex hyperbolic space;*
(3) *for any parameter $t$, $0 < t < \pi/2$, the action of the quasi-Fuchsian group $G_t$ is not quasiconformally conjugate to the action of the Fuchsian group $G = \alpha(0)$ (in both CR-structure at infinity $\partial \mathbb{H}_\mathbb{C}^2 = \mathcal{H}_2 \cup \{\infty\}$ and the complex hyperbolic space $\mathbb{H}_\mathbb{C}^2$).*

However it is still an open question whether the actions of the constructed quasi-Fuchsian groups $G_t$ and $G_{t'}$ on their limit sets $\Lambda(G_t)$ and $\Lambda(G_{t'})$ could be "quasi-conformally" conjugate, in other words, whether the canonical $G$-equivariant homeomorphism $f_{\chi_t} : \Lambda(G) \to \Lambda(G_t)$ of the limit sets (constructed in Theorem 3.2) that induces the isomorphism $\phi : G \to G_t$, $0 < t < \pi/2$, is in fact quasisymmetric.

In Sections 4 and 5, we address the basic problem of existence of non-trivial deformations of "non-real" hyperbolic manifolds (in particular, complex hyperbolic ones) and their (discrete) holonomy groups which, in contrast to the described flexibility in the real hyperbolic case, seem much more rigid. Indeed, due to Pansu [P], quasiconformal maps in the sphere at infinity of quaternionic/octanionic hyperbolic spaces induced by hyperbolic quasi-isometries are necessarily automorphisms, and thus there cannot be interesting quasiconformal deformations of corresponding structures (even any topological conjugation of two different (free) Schottky groups in those spaces cannot be quasiconformal, cf. [AM]). Secondly, due to Corlette's rigidity theorem [Co2], such manifolds are even super-rigid – analogously to Margulis super-rigidity in higher rank [M]. The last fact implies impossibility of quasi-Fuchsian deformations of such quaternionic manifolds of dimension at least 3, see [Ka]. Furthermore, complex hyperbolic manifolds share the above rigidity of quaternionic/octanionic hyperbolic manifolds. Namely, due to the Goldman's local rigidity theorem in dimension $n = 2$ [G1], every nearby discrete representation $\rho : G \to PU(2,1)$ of a cocompact lattice $G \subset PU(1,1)$ stabilizes a complex geodesic in the complex hyperbolic space $\mathbb{H}^2_\mathbb{C}$ (which is also true for small deformations of cocompact lattices $G \subset PU(n-1,1)$ in higher dimensions $n \geq 3$ [GM]), and thus the limit set $\Lambda(\rho G) \subset \partial \mathbb{H}^n_\mathbb{C}$ is always a round sphere $S^{2n-3}$. In higher dimensions $n \geq 3$, this local rigidity of complex hyperbolic $n$-manifolds $M$ homotopy equivalent to their closed complex totally geodesic hypersurfaces is even global due to a recent Yue's rigidity theorem [Y1].

Our goal here is to show that, in contrast to rigidity of complex hyperbolic n-manifolds $M$ from the above class, the Stein spaces from the following two classes of complex hyperbolic surfaces $M$ are not rigid (it seems to us that the property of a complex hyperbolic manifold to be a Stein space is crucial for its flexibility). Such a flexibility has two independent aspects related to both conditions in the Goldman's local rigidity theorem, firstly, the existence of a complex analytic subspace homotopy equivalent to the manifold $M$ and, secondly, compactness of that subspace.

Namely, as it follows from the above Theorem 4.1 and Corollary 4.2, the first class of non-rigid complex hyperbolic manifolds consists of complex Stein surfaces $M$ homotopy equivalent to their *non-compact* complex analytic subspaces (Riemann surfaces of genus $p \geq 0$ with finite hyperbolic area, with at least four punctures).

The second class of non-rigid manifolds consists of Stein spaces represented by complex hyperbolic manifolds $M$ homotopy equivalent to their closed totally *real* geodesic submanifolds. Namely, in complex dimension two, we provide a canonical construction of continuous non-trivial quasi-Fuchsian deformations of complex surfaces fibered over closed Riemannian surfaces of genus $g > 1$ depending on $3(g-1)$ continuous parameters (in addition to "Fuchsian" deformations, where in particular, the Teichmuller space of the base surface has dimension $6(g-1)$). This is the first such non-trivial deformations of fibrations with compact base (for non-compact base, see a different Goldman-Parker' deformation [GP] of ideal triangle groups $G \subset PO(2,1)$). The obtained flexibility of such holomorphic fibrations and the number of its parameters (at least $9(g-1)$) provide the first advance toward a conjecture on dimension $16(g-1)$ of the Teichmuller space of such complex surfaces. It is related to A.Weil's theorem [W] (see also [G3, p.43]), that the variety of conjugacy classes of all (not necessarily discrete) representations $G \to PU(2,1)$ near the embedding $G \subset PO(2,1)$ is a real-analytic manifold of dimension $16(g-1)$. We remark that discreteness of representations of $G \cong \pi_1 M$ is an essential condition for deformation of a complex manifold $M$ which does not follow from the mentioned Weil's result.

Our construction here is inspired by the well know bending deformations of real hyperbolic (conformal) manifolds along totally geodesic hypersurfaces. In the case of complex hyperbolic (and Cauchy-Riemannian) structures, it works however in a different way than that in the real hyperbolic case. Namely our complex bending deformations involve simultaneous bending of the base of the fibration of the complex surface $M$ as well as bendings of each of its totally geodesic fibers (see Remark 5.4). Such bending deformations of complex surfaces are associated to their real simple closed geodesics (of real codimension 3), but have nothing common with cone deformations of real hyperbolic 3-manifolds along closed geodesics (see [A4, A5]).

Furthermore, there are well known complications (cf. [KR3, P, Cap]) in constructing equivariant quasiconformal homeomorphisms in the complex hyperbolic space and in Cauchy-Riemannian geometry, which are due to necessary conditions for such maps to preserve the Kähler and contact structures (correspondingly in the complex hyperbolic space and at its infinity, the one-point compactification of the Heisenberg group$\mathcal{H}_n$). Despite that, as it follows from our construction, the complex bending deformations are induced by equivariant homeomorphisms which are in addition quasiconformal with respect to the corresponding metrics. One of our main results in this direction may be formulated as follows.

**Theorem 5.1.** *Let $G \subset PO(2,1) \subset PU(2,1)$ be a given (non-elementary) discrete group. Then, for any simple closed geodesic $\alpha$ in the Riemann 2-surface $S = H_\mathbb{R}^2/G$ and a sufficiently small $\eta_0 > 0$, there is a holomorphic family of $G$-equivariant quasiconformal homeomorphisms $F_\eta : \overline{\mathbb{H}_\mathbb{C}^2} \to \overline{\mathbb{H}_\mathbb{C}^2}$, $-\eta_0 < \eta < \eta_0$, which defines a bending (quasi-Fuchsian) deformation $\mathcal{B}_\alpha : (-\eta_0, \eta_0) \to \mathcal{R}_0(G)$ of the group $G$ along the geodesic $\alpha$, $\mathcal{B}_\alpha(\eta) = F_\eta^*$.*

We notice that such complex bending deformations depend on many independent parameters, as it is shown by application of our construction and Élie Cartan [Car] angular invariant in Cauchy-Riemannian geometry:

**Corollary 5.2.** *Let $S_p = \mathbb{H}_\mathbb{R}^2/G$ be a closed totally real geodesic surface of genus $p > 1$ in a given complex hyperbolic surface $M = \mathbb{H}_\mathbb{C}^2/G$, $G \subset PO(2,1) \subset PU(2,1)$. Then there is a real analytic embedding $\pi \circ \mathcal{B} : B^{3p-3} \hookrightarrow \mathcal{T}(M)$ of a real $(3p-3)$-ball into the Teichmüller space of $M$, defined by bending deformations along disjoint closed geodesics in $M$ and by the projection $\pi : \mathcal{D}(M) \to \mathcal{T}(M) = \mathcal{D}(M)/PU(2,1)$ in the development space $\mathcal{D}(M)$.*

The above embedding and the fact that the Teichmuller space of the base surface $S_p$ (totally geodesically) embedded in the complex surface $M$ is a complex manifold of dimension $3(g-1)$ show that we have in fact a real analytic embedding $B^{9p-9} \hookrightarrow \mathcal{T}(M)$ of a real $9(p-1)$-ball into the Teichmüller space of the complex hyperbolic surface $M$.

In our subsequent work, we apply the constructed bending deformations to answer a well known question about cusp groups on the boundary of the Teichmüller space of $\mathcal{T}(M)$ of a Stein complex surface $M$ fibering over a compact Riemann surface of genus $p > 1$:

**Theorem 5.6.** *Let $G \subset PO(2,1) \subset PU(2,1)$ be a uniform lattice isomorphic to the fundamental group of a closed surface $S_p$ of genus $p \geq 2$. Then, for any simple closed geodesic $\alpha \subset S_p = H^2_\mathbb{R}/G$, there is a continuous deformation $\rho_t = f_t^*$ induced by $G$-equivariant quasiconformal homeomorphisms $f_t : \overline{\mathbb{H}^2_\mathbb{C}} \to \overline{\mathbb{H}^2_\mathbb{C}}$ whose limit representation $\rho_\infty$ corresponds to a boundary cusp point of the Teichmüller space $\mathcal{T}(G)$, that is the boundary group $\rho_\infty(G)$ has an accidental parabolic element $\rho_\infty(g_\alpha)$ where $g_\alpha \in G$ represents the geodesic $\alpha \subset S_p$.*

We note that, due to our construction of such continuous quasiconformal deformations, they are independent if the corresponding geodesics $\alpha_i \subset S_p$ are disjoint. It implies the existence of a boundary group in $\partial \mathcal{T}(G)$ with "maximal" number of non-conjugate accidental parabolic subgroups:

**Corollary 5.7.** *Let $G \subset PO(2,1) \subset PU(2,1)$ be a uniform lattice isomorphic to the fundamental group of a closed surface $S_p$ of genus $p \geq 2$. Then there is a continuous deformation $R : \mathbb{R}^{3p-3} \to \mathcal{T}(G)$ whose boundary group $G_\infty = R(\infty)(G)$ has $3p-3$ non-conjugate accidental parabolic subgroups.*

*Acknowledgements.* Some parts of the paper were written during the author's stay at the Federal Universidade de Minas Gerais at Belo Horizonte/Brazil, the Mathematical Sciences Research Institite at Berkeley/CA, and the University of Tokyo. The author thanks them for the hospitality and gratefully acknowledges partial support by the National Science Foundation, the Federal Universidad de Minas Gerais and the Japan Society for the Promotion of Science. The author would like to thank Bill Goldman for many useful conversations. We owe our special thanks to Mario Carneiro and Nikolay Gusevskii for several stimulating discussions without which this paper would never be possible.

## 2. Complex hyperbolic geometry and geometrical finiteness

Here we recall some known facts (see, for example, [AX1, GP1, G4, KR1]) concerning the Kähler geometry of the complex hyperbolic space $\mathbb{H}^n_\mathbb{C}$, its link with the nilpotent geometry of the Heisenberg group $\mathcal{H}_n$ induced on each horosphere in $\mathbb{H}^n_\mathbb{C}$, and the Cauchy-Riemannian geometry (and contact structure) in the $(2n-1)$-sphere at infinity $\partial \mathbb{H}^n_\mathbb{C}$ which can be identified with the one-point compactification $\overline{\mathcal{H}}_n = \mathcal{H}_n \cup \{\infty\}$ of the Heisenberg group.

One can realize the complex hyperbolic space,

$$\mathbb{H}^n_\mathbb{C} = \{[z] \in \mathbb{CP}^n \; : \; z \in \mathbb{C}^{n,1} \,, \; \langle z, z \rangle < 0 \},$$

as the set of negative lines in the Hermitian vector space $\mathbb{C}^{n,1}$, with Hermitian structure given by the indefinite $(n,1)$-form $\langle z, w \rangle = z_1 \overline{w}_1 + \cdots + z_n \overline{w}_n - z_{n+1} \overline{w}_{n+1}$. Its boundary $\partial \mathbb{H}^n_\mathbb{C} = \{[z] \in \mathbb{CP}^{n,1} \; : \; \langle z, z \rangle = 0\}$ consists of all null lines in $\mathbb{CP}^n$ and is homeomorphic to the $(2n-1)$-sphere $S^{2n-1}$.

There are two common models of complex hyperbolic space $\mathbb{H}^n_\mathbb{C}$ as domains in $\mathbb{C}^n$, the unit ball $\mathbb{B}^n_\mathbb{C}$ and the Siegel domain $\mathfrak{S}_n$. They arise from two affine patches in the projective space $\mathbb{CP}^n$ related to $\mathbb{H}^n_\mathbb{C}$ and its boundary. Namely, embedding $\mathbb{C}^n$ onto the affine patch of $\mathbb{CP}^{n,1}$ defined by $z_{n+1} \neq 0$ (in homogeneous coordinates) as $A : \mathbb{C}^n \to \mathbb{CP}^n, z \mapsto [(z,1)]$, we may identify the unit ball $\mathbb{B}^n_\mathbb{C}(0,1) \subset \mathbb{C}^n$ with $\mathbb{H}^n_\mathbb{C} =$ $A(\mathbb{B}_\mathbb{C}^n)$. Here the metric in $\mathbb{C}^n$ is defined by the standard Hermitian form $\langle\langle,\rangle\rangle$, and the induced metric on $\mathbb{B}_\mathbb{C}^n$ is the Bergman metric (with constant holomorphic curvature -1) whose sectional curvature is between -1 and -1/4.

The full group $\operatorname{Isom}\mathbb{H}_\mathbb{C}^n$ of isometries of (the ball model of) $\mathbb{H}_\mathbb{C}^n$ is generated by the group of holomorphic automorphisms of the ball $\mathbb{B}_\mathbb{C}^n$ (=the projective unitary group $PU(n,1)$ defined by the group $U(n,1)$ of unitary automorphisms of $\mathbb{C}^{n,1}$ that preserve $\mathbb{H}_\mathbb{C}^n$), together with the antiholomorphic automorphism of $\mathbb{H}_\mathbb{C}^n$ defined by the $\mathbb{C}$-antilinear unitary automorphism of $\mathbb{C}^{n,1}$ given by the complex conjugation $z \mapsto \bar z$. The group $PU(n,1)$ can be embedded in a linear group due to A.Borel [Bor] (cf. [AX1, L.2.1]), hence any finitely generated group $G \subset PU(n,1)$ is residually finite and has a finite index torsion free subgroup. Elements $g \in PU(n,1)$ are of the following three types. If $g$ fixes a point in $\mathbb{H}_\mathbb{C}^n$, it is called *elliptic*. If $g$ has exactly one fixed point in the closure $\overline{\mathbb{H}_\mathbb{C}^n} \cong \overline{\mathbb{B}_\mathbb{C}^n}$, and it lies in $\partial\mathbb{H}_\mathbb{C}^n$, $g$ is called *parabolic*. If $g$ has exactly two fixed points, and they lie in $\partial\mathbb{H}_\mathbb{C}^n$, $g$ is called *loxodromic*. These three types exhaust all the possibilities.

The second model of $\mathbb{H}_\mathbb{C}^n$, as the Siegel domain, arises from the affine patch complimentary to a projective hyperplane $H_\infty$ which is tangent to $\partial\mathbb{H}_\mathbb{C}^n$ at a point $\infty \in \partial\mathbb{H}_\mathbb{C}^n$. For example, taking that point $\infty$ as $(0', -1, 1)$ with $0' \in \mathbb{C}^{n-1}$ and $H_\infty = \{[z] \in \mathbb{CP}^n : z_n + z_{n+1} = 0\}$, one has the map $\mathbf{S} : \mathbb{C}^n \to \mathbb{CP}^n \backslash H_\infty$ such that

$$\begin{pmatrix} z' \\ z_n \end{pmatrix} \longmapsto \begin{bmatrix} z' \\ \frac{1}{2} - z_n \\ \frac{1}{2} + z_n \end{bmatrix} \quad \text{where} \quad z' = \begin{pmatrix} z_1 \\ \vdots \\ z_{n-1} \end{pmatrix} \in \mathbb{C}^{n-1}.$$

In the obtained affine coordinates, the complex hyperbolic space is identified with the *Siegel domain*

$$\mathfrak{S}_n = \mathbf{S}^{-1}(\mathbb{H}_\mathbb{C}^n) = \{z \in \mathbb{C}^n : z_n + \overline{z}_n > \langle\langle z', z' \rangle\rangle\},$$

where the Hermitian form is $\langle \mathbf{S}(z), \mathbf{S}(w) \rangle = \langle\langle z', w' \rangle\rangle - z_n - \overline{w}_n$. The automorphism group of this affine model of $\mathbb{H}_\mathbb{C}^n$ is the group of affine transformations of $\mathbb{C}^n$ preserving $\mathfrak{S}_n$. Its stabilizer of the point $\infty$ is $\mathcal{H}_n \rtimes U(n-1) \cdot \exp(t)$ where $\mathcal{H}_n$ is its unipotent radical, the *Heisenberg group* that consists of all *Heisenberg translations*

$$T_{\xi,v} : (w', w_n) \mapsto \left( w' + \xi, w_n + \langle\langle \xi, w' \rangle\rangle + \frac{1}{2}(\langle\langle \xi, \xi \rangle\rangle - iv) \right),$$

where $w', \xi \in \mathbb{C}^{n-1}$ and $v \in \mathbb{R}$.

In particular, $\mathcal{H}_n$ acts simply transitively on $\partial\mathfrak{S}_n \backslash \{\infty\}$ and on each horosphere $H_t$ (in the complex hyperbolic space) centered at $\infty$, which in fact has the form:

$$H_t = \{(z', z_n) \in \mathfrak{S}_n : z_n + \overline{z}_n - \langle\langle z', z' \rangle\rangle = t\}, \quad t > 0.$$

On the base of that, one obtains the *upper half space model* for the complex hyperbolic space $\mathbb{H}_\mathbb{C}^n$ by identifying $\mathbb{C}^{n-1} \times \mathbb{R} \times [0, \infty)$ and $\overline{\mathbb{H}_\mathbb{C}^n} \backslash \{\infty\}$ as

$$(\xi, v, u) \longmapsto \begin{bmatrix} \xi \\ \frac{1}{2}(1 - \langle\langle \xi, \xi \rangle\rangle - u + iv) \\ \frac{1}{2}(1 + \langle\langle \xi, \xi \rangle\rangle + u - iv) \end{bmatrix} \in \partial\mathfrak{S}_n \backslash \{\infty\},$$

where $(\xi, v, u) \in \mathbb{C}^{n-1} \times \mathbb{R} \times [0, \infty)$ are the horospherical coordinates of the corresponding point in $\overline{\mathbb{H}_{\mathbb{C}}^n} \backslash \{\infty\}$ (with respect to the point $\infty \in \partial \mathbb{H}_{\mathbb{C}}^n$, see [GP1]).

We notice that, under this identification, the geodesics running to $\infty$ are the vertical lines $c_{\xi, v}(t) = (\xi, v, e^{2t})$ passing through points $(\xi, v) \in \mathbb{C}^{n-1} \times \mathbb{R}$. Also we see that, via the geodesic perspective from $\infty$, the "boundary plane" $H_0 = \mathbb{C}^{n-1} \times \mathbb{R} \times \{0\} = \partial \mathbb{H}_{\mathbb{C}}^n \backslash \{\infty\}$ and various horospheres correspond as $H_t \to H_u$ with $(\xi, v, t) \mapsto (\xi, v, u)$. Each of them can be identified with the Heisenberg group $\mathcal{H}_n = \mathbb{C}^{n-1} \times \mathbb{R}$. It is a 2-step nilpotent group with center $\{0\} \times \mathbb{R} \subset \mathbb{C}^{n-1} \times \mathbb{R}$, with the isometric action on itself and on $\mathbb{H}_{\mathbb{C}}^n$ by left translations:

$$T_{(\xi_0, v_0)} : (\xi, v, u) \longmapsto (\xi_0 + \xi, v_0 + v + 2 \operatorname{Im} \langle\langle \xi_0, \xi \rangle\rangle, u),$$

and the inverse of $(\xi, v)$ is $(\xi, v)^{-1} = (-\xi, -v)$. The unitary group $U(n-1)$ acts on $\mathcal{H}_n$ and $\mathbb{H}_{\mathbb{C}}^n$ by rotations: $A(\xi, v, u) = (A\xi, v, u)$ for $A \in U(n-1)$. The semidirect product $\mathcal{H}(n) = \mathcal{H}_n \rtimes U(n-1)$ is naturally embedded in $U(n, 1)$ as follows:

$$A \longmapsto \begin{pmatrix} A & 0 & 0 \\ 0 & 1 & 0 \\ 0 & 0 & 1 \end{pmatrix} \in U(n, 1) \quad \text{for} \quad A \in U(n-1),$$

$$(\xi, v) \longmapsto \begin{pmatrix} I_{n-1} & \xi & \xi \\ -\bar{\xi}^t & 1 - \frac{1}{2}(|\xi|^2 - iv) & -\frac{1}{2}(|\xi|^2 - iv) \\ \bar{\xi}^t & \frac{1}{2}(|\xi|^2 - iv) & 1 + \frac{1}{2}(|\xi|^2 - iv) \end{pmatrix} \in U(n, 1)$$

where $(\xi, v) \in \mathcal{H}_n = \mathbb{C}^{n-1} \times \mathbb{R}$ and $\bar{\xi}^t$ is the conjugate transpose of $\xi$.

The action of $\mathcal{H}(n)$ on $\overline{\mathbb{H}_{\mathbb{C}}^n} \backslash \{\infty\}$ also preserves the Cygan metric $\rho_c$ there, which plays the same role as the Euclidean metric does on the upper half-space model of the real hyperbolic space $\mathbb{H}^n = \mathbb{H}_{\mathbb{R}}^n$ and is induced by the following norm:

$$||(\xi, v, u)||_c = |\, ||\xi||^2 + u - iv|^{1/2}, \quad (\xi, v, u) \in \mathbb{C}^{n-1} \times \mathbb{R} \times [0, \infty). \tag{2.1}$$

The relevant geometry on each horosphere $H_u \subset \mathbb{H}_{\mathbb{C}}^n$, $H_u \cong \mathcal{H}_n = \mathbb{C}^{n-1} \times \mathbb{R}$, is the spherical $CR$-geometry induced by the complex hyperbolic structure. The geodesic perspective from $\infty$ defines $CR$-maps between horospheres, which extend to $CR$-maps between the one-point compactifications $H_u \cup \infty \approx S^{2n-1}$. In the limit, the induced metrics on horospheres fail to converge but the $CR$-structure remains fixed. In this way, the complex hyperbolic geometry induces $CR$-geometry on the sphere at infinity $\partial \mathbb{H}_{\mathbb{C}}^n \approx S^{2n-1}$, naturally identified with the one-point compactification of the Heisenberg group $\mathcal{H}_n$.

Our main assumption on a complex hyperbolic $n$-manifold $M$ is the geometrical finiteness of its fundamental group $\pi_1(M) = G \subset \operatorname{Isom} \mathbb{H}_{\mathbb{C}}^n$, which in particular implies that the discrete group $G$ is finitely presented [AX1].

Here a subgroup $G \subset \operatorname{Isom} \mathbb{H}_{\mathbb{C}}^n$ is called *discrete* if it is a discrete subset of $\operatorname{Isom} \mathbb{H}_{\mathbb{C}}^n$. The *limit set* $\Lambda(G) \subset \partial \mathbb{H}_{\mathbb{C}}^n$ of a discrete group $G$ is the set of accumulation points of (any) orbit $G(y)$, $y \in \mathbb{H}_{\mathbb{C}}^n$. The complement of $\Lambda(G)$ in $\partial \mathbb{H}_{\mathbb{C}}^n$ is called the *discontinuity set* $\Omega(G)$. A discrete group $G$ is called *elementary* if its limit set $\Lambda(G)$ consists of at most two points. An infinite discrete group $G$ is called *parabolic* if it has exactly one fixed point $\operatorname{fix}(G)$; then $\Lambda(G) = \operatorname{fix}(G)$, and $G$ consists of either parabolic or elliptic elements. As it was observed by many authors, parabolicity in the variable curvature case is not as easy a condition to deal with as it is in the constant curvature space. Even the notion of a parabolic cusp point become somewhat complicated. Namely, following to [Bow], a parabolic fixed point $p$ of a discrete group $G \subset \text{Isom}\,\mathbb{H}_{\mathbb{C}}^n$ is called a *cusp point* if the quotient $(\Lambda(G)\setminus\{p\})/G_p$ of the limit set of $G$ by the action of the parabolic stabilizer $G_p = \{g \in G : g(p) = p\}$ is compact. However our approach [AX1-AX3] makes this notion and geometrical finiteness in pinched negative curvature itself much more transparent.

Geometrical finiteness has been essentially used for real hyperbolic manifolds, where geometric analysis and ideas of Thurston provided powerful tools for understanding of their structure. Due to the absence of totally geodesic hypersurfaces in a space of variable negative curvature and recent results [AX1, GP1] on Dirichlet polyhedra for simplest parabolic groups in $\mathbb{H}_{\mathbb{C}}^n$, we cannot use the original definition of geometrical finiteness which came from an assumption that the corresponding real hyperbolic manifold $M = \mathbb{H}^n/G$ may be decomposed into a cell by cutting along a finite number of its totally geodesic hypersurfaces, that is the group $G$ should possess a finite-sided fundamental polyhedron, see [Ah]. However, we can use many other (equivalent) definitions of geometrical finiteness.

The first one, **GF1** (originally due to A.Beardon and B.Maskit [BM]) defines *geometrically finite* discrete groups $G \subset \text{Isom}\,\mathbb{H}_{\mathbb{C}}^n$ (and their complex hyperbolic orbifolds $M = \mathbb{H}_{\mathbb{C}}^n/G$) as those whose limit set $\Lambda(G)$ entirely consists of parabolic cusps and conical limit points. Here a limit point $z \in \Lambda(G)$ is called a *conical limit point* of $G$ if, for some (and hence every) geodesic ray $\ell \subset \mathbb{H}_{\mathbb{C}}^n$ ending at $z$, there is a compact set $K \subset \mathbb{H}_{\mathbb{C}}^n$ such that $g(\ell) \cap K \neq \emptyset$ for infinitely many elements $g \in G$. The last condition is equivalent to the following [BM, AX3]:

For every geodesic ray $\ell \subset \mathbb{H}_{\mathbb{C}}^n$ ending at $z$ and for every $\delta > 0$, there is a point $x \in \mathbb{H}_{\mathbb{C}}^n$ and a sequence of distinct elements $g_i \in G$ such that the orbit $\{g_i(x)\}$ approximates $z$ inside the $\delta$-neighborhood $N_\delta(\ell)$ of the ray $\ell$.

Our study of geometrical finiteness in variable curvature [AX1-AX3] is based on analysis of geometry and topology of thin (parabolic) ends of corresponding manifolds and parabolic cusps of discrete isometry groups $G \subset PU(n,1)$. Namely, suppose a point $p \in \partial \mathbb{H}_{\mathbb{C}}^n$ is fixed by some parabolic element of a given discrete group $G \subset \text{Isom}\,\mathbb{H}_{\mathbb{C}}^n$, and $G_p$ is the stabilizer of $p$ in $G$. Conjugating $G$ by an element $h_p \in PU(n,1), h_p(p) = \infty$, we may assume that the stabilizer $G_p$ is a subgroup $G_\infty \subset \mathcal{H}(n) = \mathcal{H}_n \rtimes U(n-1)$. In particular, if $p$ is the origin $0 \in \mathcal{H}_n$, the transformation $h_0$ can be taken as the Heisenberg inversion $\mathcal{I}$ in the hyperchain $\partial \mathbb{H}_{\mathbb{C}}^{n-1}$. It preserves the unit Heisenberg sphere $S_c(0,1) = \{(\xi,v) \in \mathcal{H}_n : \|(\xi,v)\|_c = 1\}$ and acts in $\mathcal{H}_n$ as follows:

$$\mathcal{I}(\xi,v) = \left(\frac{\xi}{|\xi|^2 - iv}, \frac{-v}{v^2 + |\xi|^4}\right) \quad \text{where } (\xi,v) \in \mathcal{H}_n = \mathbb{C}^{n-1} \times \mathbb{R}. \tag{2.1}$$

For any other point $p$, we may take $h_p$ as the Heisenberg inversion $\mathcal{I}_p$ which preserves the unit Heisenberg sphere $S_c(p,1) = \{(\xi,v) : \rho_c(p,(\xi,v)) = 1\}$ centered at $p$. The inversion $\mathcal{I}_p$ is the conjugate of $\mathcal{I}$ by the Heisenberg translation $T_p$; it maps $p$ to $\infty$.

After such a conjugation, we can regard the parabolic stabilizer $G_p$ as a discrete isometry group acting in the (nilpotent) Heisenberg group $\mathcal{H}_n$. This action is completely described by our following result [AX1-AX3]:

**Theorem 2.1.** *Let $\mathcal{N}$ be a connected, simply connected nilpotent Lie group, $C$ be a compact group of automorphisms of $\mathcal{N}$, and $\Gamma \subset \mathcal{N} \rtimes C$ be a discrete subgroup.*

*Then there exist a connected Lie subgroup $\mathcal{N}_\Gamma$ of $\mathcal{N}$ and a finite index subgroup $\Gamma^*$ of $\Gamma$ with the following properties:*

(1) *There exists $b \in \mathcal{N}$ such that $b\Gamma b^{-1}$ preserves $\mathcal{N}_\Gamma$;*
(2) $\mathcal{N}_\Gamma / b\Gamma b^{-1}$ *is compact;*
(3) $b\Gamma^* b^{-1}$ *acts on $\mathcal{N}_\Gamma$ by left translations and the action of $b\Gamma^* b^{-1}$ on $\mathcal{N}_\Gamma$ is free.*

Due to this Theorem, there is a connected Lie subgroup $\mathcal{H}_\infty \subseteq \mathcal{H}_n$ preserved by $G_\infty$ (up to changing the origin). So we can make the following definition.

**Definition 2.2.** A set $U_{p,r} \subset \overline{\mathbb{H}_{\mathbb{C}}^n} \setminus \{p\}$ is called a *standard cusp neighborhood of radius $r > 0$* at a parabolic fixed point $p \in \partial \mathbb{H}_{\mathbb{C}}^n$ of a discrete group $G \subset PU(n,1)$ if, for the Heisenberg inversion $\mathcal{I}_p \in PU(n,1)$ with respect to the unit sphere $S_c(p,1)$, $\mathcal{I}_p(p) = \infty$, the following conditions hold:

(1) $U_{p,r} = \mathcal{I}_p^{-1}\left(\{x \in \mathbb{H}_{\mathbb{C}}^n \cup \mathcal{H}_n \ : \ \rho_c(x, \mathcal{H}_\infty) \geq 1/r\}\right)$;
(2) $U_{p,r}$ is precisely invariant with respect to $G_p \subset G$, that is:

$$\gamma(U_{p,r}) = U_{p,r} \quad \text{for} \quad \gamma \in G_p \quad \text{and} \quad g(U_{p,r}) \cap U_{p,r} = \emptyset \quad \text{for} \quad g \in G \backslash G_p.$$

Now, due to [AX1], we can give a geometric definition of a cusp point. Namely, a parabolic fixed point $p \in \partial \mathbb{H}_{\mathbb{C}}^n$ of a discrete group $G \subset \text{Isom}\,\mathbb{H}_{\mathbb{C}}^n$ is a cusp point if and only if it has a standard cusp neighborhood $U_{p,r} \subset \overline{\mathbb{H}_{\mathbb{C}}^n} \setminus \{p\}$.

This fact and [Bow] allow us to give another equivalent definitions of geometrical finiteness which is originally due to A.Marden [Ma]. In particular it follows that a discrete subgroup $G$ in $PU(n,1)$ is *geometrically finite* (**GF2**) if and only if its quotient space

$$M(G) = [\mathbb{H}_{\mathbb{C}}^n \cup \Omega(G)]/G \tag{2.2}$$

has finitely many ends, and each of them is a cusp end, that is an end whose neighborhoods can be taken (for an appropriate $r > 0$) in the form:

$$U_{p,r}/G_p \approx (S_{p,r}/G_p) \times (0,1], \tag{2.3}$$

where

$$S_{p,r} = \partial_H U_{p,r} = \mathcal{I}_p^{-1}\left(\{x \in H_{\mathbb{C}}^n \cup \mathcal{H}_n \ : \ \rho_c(x, \mathcal{H}_\infty) = 1/r\}\right).$$

Now we see that a geometrically finite manifold can be decomposed into a compact submanifold and finitely many cusp submanifolds of the form (2.3). Clearly, each of such cusp ends is homotopy equivalent to a Heisenberg $(2n-1)$-manifold which can be described as follows [AX1]:

**Theorem 2.3.** *Let $\Gamma \subset \mathcal{H}_n \rtimes U(n-1)$ be a torsion-free discrete group acting on the Heisenberg group $\mathcal{H}_n = \mathbb{C}^{n-1} \times \mathbb{R}$ with non-compact quotient. Then the quotient $\mathcal{H}_n/\Gamma$ has zero Euler characteristic and is a vector bundle over a compact manifold. Furthermore, this compact manifold is finitely covered by a nil-manifold which is either a torus or the total space of a circle bundle over a torus.*

Now it follows that each cusp end is homotopy equivalent to a compact $k$-manifold, $k \leq 2n - 1$, finitely covered by a nil-manifold which is either a (flat) torus or the total space of a circle bundle over a torus. It implies that the fundamental groups of cusp ends are finitely presented, and we get the following finiteness [AX1]:

**Corollary 2.4.** *Geometrically finite groups $G \subset \mathrm{Isom}\,\mathbb{H}_\mathbb{C}^n$ are finitely presented.*

Another application of our geometric approach shows that cusp ends of a geometrically finite complex hyperbolic orbifolds $M$ have, up to a finite covering of $M$, a very simple structure [AX1]:

**Theorem 2.5.** *Let $G \subset \mathrm{Isom}\,\mathbb{H}_\mathbb{C}^n$ be a geometrically finite discrete group. Then $G$ has a subgroup $G_0$ of finite index such that every parabolic subgroup of $G_0$ is isomorphic to a discrete subgroup of the Heisenberg group $\mathcal{H}_n = \mathbb{C}^{n-1} \times \mathbb{R}$. In particular, each parabolic subgroup of $G_0$ is free Abelian or 2-step nilpotent.*

In terms of finite coverings, this result has the following sense:

**Corollary 2.6.** *For a given geometrically finite orbifold $M(G) = \overline{\mathbb{H}_\mathbb{C}^n}\backslash\Lambda(G)/G$, there is its finite covering $\hat{M}$ such that neighborhoods of each (cusp) end of $\hat{M}$ are homeomorphic to the product of infinite interval $[0,\infty)$, a closed k-dimensional ball $B^k$ and a closed $(2n-k-1)$-dimensional manifold which is either the (flat) torus $T^{2n-k-1}$ or the total space of a (non-trivial) circle bundle over the torus $T^{2n-k-2}$.*

Finally we formulate two additional definitions of geometrical finiteness which are originally due to W.Thurston [T]:

**(GF3):** The thick part of the minimal convex retract (=convex core) $C(G)$ of $\mathbb{H}_\mathbb{C}^n/G$ is compact.

**(GF4):** For some $\epsilon > 0$, the uniform $\epsilon$-neighborhood of the convex core $C(G) \subset \mathbb{H}_\mathbb{C}^n/G$ has finite volume, and there is a bound on the orders of finite subgroups of $G$.

Due to Bowditch [Bow], the four definitions **GF1, GF2, GF3** and **GF4** of geometrical finiteness for a discrete group $G \subset \mathrm{Isom}\,\mathbb{H}_\mathbb{C}^n$ are all equivalent, see also [AX1, AX3].

Now we would like to define the above terms of "convex core" and "thick part" of a complex hyperbolic orbifold $M$. Namely, the convex core $C(G)$ of a complex hyperbolic orbifold $M = \mathbb{H}_\mathbb{C}^n/G$ can be obtained as the $G$-quotient of the complex hyperbolic convex hull $C(\Lambda(G))$ of the limit set $\Lambda(G)$. Here the convex hull $C(\Lambda(G)) \subset \mathbb{H}_\mathbb{C}^n$ is the minimal convex subset in $\mathbb{H}_\mathbb{C}^n$ whose closure in $\overline{\mathbb{H}_\mathbb{C}^n}$ contains the limit set $\Lambda(G)$. Clearly, it is $G$-invariant, and its quotient $C(G) = C(\Lambda(G))/G$ is the minimal convex retract of $\mathbb{H}_\mathbb{C}^n/G$; we call it the *convex core* of $M = \mathbb{H}_\mathbb{C}^n/G$.

Now let $\epsilon$ be any positive number less than the Margulis constant in dimension $n$, $\epsilon(n)$. Then for a given discrete group $G \subset \mathrm{Isom}\,\mathbb{H}_\mathbb{C}^n$ and its orbifold $M = \mathbb{H}_\mathbb{C}^n/G$, the $\epsilon$-thin part $\mathrm{thin}_\epsilon(M)$ is defined as:

$$\mathrm{thin}_\epsilon(M) = \{x \in \mathbb{H}_\mathbb{C}^n : G_\epsilon(x) = \langle g \in G : d(x,g(x)) < \epsilon \rangle \text{ is infinite}\}/G. \quad (2.4)$$

The $\epsilon$-thick part $\mathrm{thick}_\epsilon(M)$ of an orbifold $M$ is defined as the closure of the complement to the $\epsilon$-thin part $\mathrm{thin}_\epsilon(M) \subset M$.

As a consequence of the Margulis Lemma [M, BGS], there is the following description [BGS, Bow] of the thin part of a negatively curved orbifold which we formulate for complex hyperbolic geometry:

**Theorem 2.7.** *Let $G \subset \mathrm{Isom}\,\mathbb{H}_\mathbb{C}^n$ be a discrete group and $\epsilon$, $0 < \epsilon < \epsilon(n)$, be choosen. Then the $\epsilon$-thin part $\mathrm{thin}_\epsilon(M)$ of $M = \mathbb{H}_\mathbb{C}^n/G$ is a disjoint union of its connected components, and each such component has the form $T_\epsilon(\Gamma)/\Gamma$ where $\Gamma$* is a maximal infinite elementary subgroup of $G$. Here, for each such elementary subgroup $\Gamma \subset G$, the connected component (Margulis region)

$$T_\epsilon = \{x \in \mathbb{H}^n_{\mathbb{C}} : \Gamma_\epsilon(x) = \langle g \in \Gamma : d(x, \gamma(x)) < \epsilon \rangle \text{ is infinite}\} \quad (2.5)$$

is precisely invariant for $\Gamma$ in $G$:

$$\Gamma(T_\epsilon) = T_\epsilon, \quad g(T_\epsilon) \cap T_\epsilon = \emptyset \quad \text{for any } g \in G \backslash \Gamma. \quad (2.6)$$

We note that in the real hyperbolic case of dimension 2 and 3, a Margulis region $T_\epsilon$ in (2.5) with parabolic stabilizer $\Gamma \subset G$ can be taken as a horoball neighborhood centered at the parabolic fixed point $p$, $\Gamma(p) = p$. It is not true in general because of Apanasov's construction [A3] in real hyperbolic spaces of dimension at least 4. As we discussed it in [AX1], this construction works in complex hyperbolic spaces $\mathbb{H}^n_{\mathbb{C}}$, $n \geq 2$, as well.

However, applying the structural Theorem 2.1 to actions of parabolic groups nearby their fixed points, we obtain a description of parabolic Margulis regions for any discrete groups $G \subset \text{Isom } \mathbb{H}^n_{\mathbb{C}}$ (even in more general situation of pinched Hadamard manifolds, see [AX3, Lemma 5.2]). Namely, let $\Gamma \subset G$ be such a discrete parabolic subgroup. Without loss of generality, we may assume that its fixed point $p \in \partial \mathbb{H}^n_{\mathbb{C}}$ is $\infty$ in the Siegel domain, or equivalently, in the upper half-space model of $\mathbb{H}^n_{\mathbb{C}}$. Then, on each horosphere $H_t \subset \mathbb{H}^n_{\mathbb{C}}$ centered at $p = \infty$, the parabolic group $\Gamma$ acts as a discrete subgroup of $\mathcal{H}_n \rtimes U(n-1)$. Hence, applying Theorem 2.1, we have a $\Gamma$-invariant connected subspace $\mathcal{H}_\Gamma \subset \partial \mathbb{H}^n_{\mathbb{C}} \backslash \{p\}$ where $\Gamma$ acts co-compactly, and on which a finite index subgroup $\Gamma^* \subset \Gamma$ acts freely by left translations. We define the subspace $\tau_\Gamma \subset \mathbb{H}^n_{\mathbb{C}}$ to be spanned by $\mathcal{H}_\Gamma$ and all geodesics $(z, p) \subset \mathbb{H}^n_{\mathbb{C}}$ that connect $z \in \mathcal{H}_\Gamma$ to the parabolic fixed point $p$. Let $\tau_\Gamma^t$ be the "half-plane" in $\tau_\Gamma$ of a height $t > 0$, that is the part of $\tau_\Gamma$ whose last horospherical coordinate is at least $t$. Then, due to [AX3, Lemma 5.2], we have:

**Lemma 2.8.** *Let $G \subset \text{Isom } \mathbb{H}^n_{\mathbb{C}}$ be a discrete group and $p$ a parabolic fixed point of $G$. Let $T_\epsilon$ be a Margulis region for $p$ as given in (2.5) and let $\tau_\Gamma^t$ be the half-plane defined as above. Then for any $\delta$, $0 < \delta < \epsilon/2$, there exists a large enough number $t > 0$ such that the Margulis region $T_\epsilon$ contains the $\delta$-neighborhood $N_\delta(\tau_\Gamma^t)$ of the half-plane $\tau_\Gamma^t$.*

This fact, Theorem 2.7 and **GF3**-characterization of geometrical finiteness (compactness of the $\epsilon$-thick part of the convex core $C(G)$ with sufficiently small $\epsilon > 0$) imply the following (equivalent) description of geometrically finite complex hyperbolic orbifolds $M = \mathbb{H}^n_{\mathbb{C}}/G$. Namely it follows that the action of a geometrically finite discrete group $G \subset \text{Isom } \mathbb{H}^n_{\mathbb{C}}$ on the convex hull $C(\Lambda(G))$ has a $G$-invariant family of precisely invariant disjoint horoballs centered at parabolic fixed points (their sufficiently small sizes are determined by Lemma 2.8). In other words, we have:

**Corollary 2.9.** *A discrete group $G \subset \text{Isom } \mathbb{H}^n_{\mathbb{C}}$ is geometrically finite if and only if there is a $G$-invariant family of disjoint open horoballs $B_i \subset \mathbb{H}^n_{\mathbb{C}}$ centered at parabolic fixed points $p_i \in \partial \mathbb{H}^n_{\mathbb{C}}$ of the group $G$ such that the orbifold*

$$C_0(G) = \big(C(\Lambda(G)) \backslash \cup_i B_i\big)/G \quad (2.7)$$

*is compact (and homotopy equivalent to $M = \mathbb{H}^n_{\mathbb{C}}/G$).*

3. Geometric isomorphisms

Here we would like to discuss the well known problem of geometric realizations of isomorphisms of discrete groups. Adapting its formulation in §1 for discrete groups $G, H \subset PU(n,1)$, we have:

**Problem 3.1.** *Given a type preserving isomorphism $\varphi : G \to H$ of discrete groups $G, H \subset PU(n,1)$, find subsets $X_G, X_H \subset \overline{\mathbb{H}_{\mathbb{C}}^n}$ invariant for the action of groups $G$ and $H$, respectively, and an equivariant homeomorphism $f_\varphi : X_G \to X_H$ which induces the isomorphism $\varphi$. Determine metric properties of $f_\varphi$, in particular, whether it is either quasisymmetric or quasiconformal with respect to either the Bergman metric in $\mathbb{H}_{\mathbb{C}}^n$ or the induced Cauchy-Riemannian structure at infinity $\partial \mathbb{H}_{\mathbb{C}}^n$.*

Such type problems were studied by several authors. In the case of lattices $G$ and $H$ in rank 1 symmetric spaces $X$, G.Mostow [Mo1] proved in his celebrated rigidity theorem that such isomorphisms $\varphi : G \to H$ can be extended to inner isomorphisms of $X$, provided that there is no analytic homomorphism of $X$ onto $PSL(2,\mathbb{R})$. For that proof, it was essential to prove that $\varphi$ can be induced by a quasiconformal homeomorphism of the sphere at infinity $\partial X$ which is the one point compactification of a (nilpotent) Carnot group $N$ (for quasiconformal mappings in Heisenberg and Carnot groups, see [KR1, KR2, P]).

If geometrically finite groups $G, H \subset PU(n,1)$ have parabolic elements and are neither lattices nor trivial, the only known geometric realization of their isomorphisms in dimension $\dim X \geq 3$ is due to P.Tukia's [Tu] isomorphism theorem for real hyperbolic spaces $X = \mathbb{H}_{\mathbb{R}}^n$. However, that Tukia's construction (based on geometry of convex hulls of the limit sets $\Lambda(G)$ and $\Lambda(H)$) cannot be used in the case of variable negative curvature due to lack of control over convex hulls (where the convex hull of three points may be 4-dimensional), especially nearby parabolic fixed points. Another (dynamical) approach due to C.Yue [Y2, Cor.B] (and the Anosov-Smale stability theorem for hyperbolic flows) can be used only for convex cocompact groups $G$ and $H$ [Y3]. As a first step in solving the above geometrization Problem 3.1, we have the following isomorphism theorem [A9, A11] in the complex hyperbolic space:

**Theorem 3.2.** *Let $\phi : G \to H$ be a type preserving isomorphism of two non-elementary geometrically finite groups $G, H \subset \operatorname{Isom} \mathbb{H}_{\mathbb{C}}^n$. Then there exists a unique equivariant homeomorphism $f_\phi : \Lambda(G) \to \Lambda(H)$ of their limit sets that induces the isomorphism $\phi$. Moreover, if $\Lambda(G) = \partial \mathbb{H}_{\mathbb{C}}^n$, the homeomorphism $f_\phi$ is the restriction of a hyperbolic isometry $h \in \operatorname{Isom} \mathbb{H}_{\mathbb{C}}^n$.*

*Proof.* For completeness, we prove this fact (following to [A9, A11]). We consider the Cayley graph $K(G, \sigma)$ of a group $G$ with a given finite set $\sigma$ of generators. This is a 1-complex whose vertices are elements of $G$, and such that two vertices $a, b \in G$ are joined by an edge if and only if $a = bg^{\pm 1}$ for some generator $g \in \sigma$. Let $|*|$ be the word norm on the graph $K(G, \sigma)$, that is the norm $|g|$ equals the minimal length of words in the alphabet $\sigma$ representing a given element $g \in G$. Choosing a function $\rho$ such that $\rho(r) = 1/r^2$ for $r > 0$ and $\rho(0) = 1$, one can define the length of an edge $[a,b] \subset K(G, \sigma)$ as $d_\rho(a,b) = \min\{\rho(|a|), \rho(|b|)\}$. Considering paths of minimal length in the sense of the function $d_\rho(a,b)$, one can extend it to a metric on the Cayley graph $K(G, \sigma)$. So taking the Cauchy completion $\overline{K(G, \sigma)}$ of that metric space, we have the definition of the group completion $\overline{G}$ as the compact metric space $\overline{K(G,\sigma)}\backslash K(G,\sigma)$, see [Fl]. Up to a Lipschitz equivalence, this definition does not depend on $\sigma$. It is also clear that, for a cyclic group $\mathbb{Z}$, its completion $\overline{\mathbb{Z}}$ consists of two points. Nevertheless, for a nilpotent group $G$ with one end, its completion $\overline{G}$ is a one-point set [Fl].

Now we can define a proper equivariant embedding $F : K(G,\sigma) \hookrightarrow \mathbb{H}_{\mathbb{C}}^n$ of the Cayley graph of a given geometrically finite group $G \subset PU(n,1)$. To do that we may assume that the stabilizer of a base point, say $0 \in \mathbb{B}_{\mathbb{C}}^n \cong \mathbb{H}_{\mathbb{C}}^n$, is trivial. Then we set $F(g) = g(0)$ for any vertex $g \in K(G,\sigma)$, and $F$ maps any edge $[a,b] \subset K(G,\sigma)$ to the geodesic segment $[a(0), b(0)] \subset \mathbb{H}_{\mathbb{C}}^n$ joining the points $a(0)$ and $b(0)$.

**Proposition 3.3.** *For a geometrically finite discrete group $G \subset \mathrm{Isom}\,\mathbb{H}_{\mathbb{C}}^n$, there are constants $K, K' > 0$ such that the following bounds hold for all elements $g \in G$ with $|g| \geq K'$:*

$$\ln(2|g| - K)^2 - \ln K^2 \leq d(0, g(0)) \leq K|g|\,. \qquad (3.4)$$

T. he proof of this claim is based on a comparison of the Bergman metric $d(*,*)$ and the path metric $d_0(*,*)$ on the following subset $\mathbb{H}_0 \subset \mathbb{H}_{\mathbb{C}}^n$. As in §2, let $C(\Lambda(G)) \subset \mathbb{H}_{\mathbb{C}}^n$ be the convex hull of the limit set $\Lambda(G) \subset \partial \mathbb{H}_{\mathbb{C}}^n$ of the group $G$. Since $G$ is geometrically finite, the complement in $M(G)$ to neighbourhoods of (finitely many) cusp ends is compact, and its retract can be taken as the compact suborbifold $C_0(G)$ in the convex core $C(G)$, see (2.7) and Corollary 2.9. Its universal cover $\mathbb{H}_0 \subset C(\Lambda(G))$ is the complement in the convex hull $C(\Lambda(G))$ to a $G$-invariant family of disjoint open horoballs $B_i \subset \mathbb{H}_{\mathbb{C}}^n$ centered at parabolic fixed points $p_i \in \partial \mathbb{H}_{\mathbb{C}}^n$ of the group $G$, and each of which is precisely invariant with respect to its (parabolic) stabilizer $G_i \subset G$.

Now, having a co-compact action of the group $G$ on the domain $\mathbb{H}_0 \subset \mathbb{H}_{\mathbb{C}}^n$ whose boundary includes some horospheres, we can reduce our comparison of distance functions $d = d(x, x')$ and $d_0 = d_0(x, x')$ to their comparison on a horosphere. So we can take points $x = (0, 0, u)$ and $x' = (\xi, v, u)$ on a "horizontal" horosphere $H_u = \mathbb{C}^{n-1} \times \mathbb{R} \times \{u\} \subset \mathbb{H}_{\mathbb{C}}^n$. Then the distances $d$ and $d_0$ are as follows [Pr2]:

$$\cosh^2 \frac{d}{2} = \frac{1}{4u^2}\left(|\xi|^4 + 4u|\xi|^2 + 4u^2 + v^2\right),\quad d_0^2 = \frac{|\xi|^2}{u} + \frac{v^2}{4u^2}\,. \qquad (3.5)$$

This comparison and the basic fact due to Cannon [Can] that, for a co-compact action of a group $G$ in a metric space $X$, its Cayley graph can be quasi-isometrically embedded into $X$, finish our proof of (3.4), compare [A7].

Now we apply Proposition 3.3 to define a $G$-equivariant extension of the map $F$ from the Cayley graph $K(G,\sigma)$ to the group completion $\overline{G}$. Since the group completion of any parabolic subgroup $G_p \subset G$ is either a point or a two-point set (depending on whether $G_p$ is a finite extension of cyclic or a nilpotent group with one end), we get

**Theorem 3.4.** *For a geometrically finite discrete group $G \subset \mathrm{Isom}\,\mathbb{H}_{\mathbb{C}}^n$, there is a continuous $G$-equivariant map $\Phi_G : \overline{G} \to \Lambda(G)$. Moreover, the map $\Phi_G$ is bijective everywhere but the set of parabolic fixed points $p \in \Lambda(G)$ whose stabilizers $G_p \subset G$ have rank one. On this set, the map $\Phi_G$ is two-to-one.*

Now we can finish our proof of Theorem 3.2 by looking at the following diagram of maps:

$$\Lambda(G) \xleftarrow{\Phi_G} \overline{G} \xrightarrow{\overline{\phi}} \overline{H} \xrightarrow{\Phi_H} \Lambda(H)\,,$$

where the homeomorphism $\overline{\phi}$ is induced by the isomorphism $\phi$, and the continuous maps $\Phi_G$ and $\Phi_H$ are defined by Theorem 3.4. Namely, one can define a map $f_\phi = \Phi_H \overline{\phi} \Phi_G^{-1}$. Here the map $\Phi_G^{-1}$ is the right inverse to $\Phi_G$, which exists due to Theorem 3.4. Furthermore, the map $\Phi_G^{-1}$ is bijective everywhere but the set of parabolic fixed points $p \in \Lambda(G)$ whose stabilizers $G_p \subset G$ have rank one, where the map $\Phi_G^{-1}$ is 2-to-1. Hence the composition map $f_\phi$ is bijective and $G$-equivariant. Its uniqueness follows from its continuity and the fact that the image of the attractive fixed point of a loxodromic element $g \in G$ must be the attractive fixed point of the loxodromic element $\phi(g) \in H$ (such loxodromic fixed points are dense in the limit set, see [A1]).

The last claim of the Theorem 3.2 directly follows from the Mostow rigidity theorem [Mo1] because a geometrically finite group $G \subset \mathrm{Isom}\,\mathbb{H}_\mathbb{C}^n$ with $\Lambda(G) = \partial \mathbb{H}_\mathbb{C}^n$ is co-finite: $\mathrm{Vol}\,(\mathbb{H}_\mathbb{C}^n/G) < \infty$.

□

*Remark 3.5.* Our proof of Theorem 3.2 can be easily extended to the general situation of type preserving isomorphisms of geometrically finite discrete groups in pinched Hadamard manifolds due to recent results in [AX3]. Namely, it is possible to construct equivariant homeomorphisms $f_\phi : \Lambda(G) \to \Lambda(H)$ conjugating the actions (on the limit sets) of isomorphic geometrically finite groups $G, H \subset \mathrm{Isom}\,X$ in a (symmetric) space $X$ with pinched negative curvature $K$, $-b^2 \leq K \leq -a^2 < 0$. Actually, bounds similar to (3.4) in Prop. 3.3 (crucial for our argument) can be obtained from a result due to Heintze and Im Hof [HI, Th.4.6] which compares the geometry of horospheres $H_u \subset X$ with that in the spaces of constant curvature $-a^2$ and $-b^2$, respectively. It gives, that for all $x, y \in H_u$ and their distances $d = d(x, y)$ and $d_u = d_u(x, y)$ in the space $X$ and in the horosphere $H_u$, respectively, one has that

$$\frac{2}{a} \sinh(a \cdot d/2) \leq d_u \leq \frac{2}{b} \sinh(b \cdot d/2).$$

Upon existence of such a (canonical) homeomorphisms $f_\varphi$ that induces a given type-preserving isomorphisms $\varphi$ of discrete subgroups of $\mathrm{Isom}\,\mathbb{H}_\mathbb{C}^n$, the geometric realization Problem 3.1 can be reduced to the questions whether $f_\varphi$ is quasisymmetric with respect to the Carnot-Carathéodory (or Cygan) metric, and whether there exists its $G$-equivariant extension to a bigger set (in particular to the sphere at infinity $\partial X$ or even to the whole space $\overline{\mathbb{H}_\mathbb{C}^n}$, cf. [KR2]) inducing the isomorphism $\varphi$. For convex cocompact groups obtained by nearby representations, this may be seen as a generalization of D.Sullivan stability theorem [Su2], see also [A7]. We shall discuss that question in the next Section. Also we note that, besides the metrical (quasisymmetric) part of the geometrization Problem 3.1, there are some topological obstructions for extensions of equivariant homeomorphisms $f_\varphi, f_\varphi : \Lambda(G) \to \Lambda(H)$. It follows from the next example.

**Example 3.6.** *Let $G \subset PU(1,1) \subset PU(2,1)$ and $H \subset PO(2,1) \subset PU(2,1)$ be two geometrically finite (loxodromic) groups isomorphic to the fundamental group $\pi_1(S_g)$ of a compact oriented surface $S_g$ of genus $g > 1$. Then the equivariant homeomorphism $f_\varphi : \Lambda(G) \to \Lambda(H)$ cannot be homeomorphically extended to the whole sphere $\partial \mathbb{H}_\mathbb{C}^2 \approx S^3$.*

*Proof.* The obstruction in this example is topological and is due to the fact that the quotient manifolds $M_1 = \mathbb{H}_{\mathbb{C}}^2/G$ and $M_2 = \mathbb{H}_{\mathbb{C}}^2/H$ are not homeomorphic. Namely, these complex surfaces are disk bundles over the Riemann surface $S_g$ and have different Toledo invariants: $\tau(\mathbb{H}_{\mathbb{C}}^2/G) = 2g - 2$ and $\tau(\mathbb{H}_{\mathbb{C}}^2/H) = 0$, see [To]. □

The complex structures of the complex surfaces $M_1$ and $M_2$ are quite different, too. The first manifold $M_1$ has a natural embedding of the Riemann surface $S_g$ as a closed analytic totally geodesic submanifold, and hence $M_1$ cannot be a Stein manifold. The second manifold $M_2$, a disk bundle over the Riemann surface $S_g$ has a totally geodesic real section and is a Stein manifold due to a result by Burns–Shnider [BS].

Moreover due to Goldman [G1], since the surface $S_g \subset M_1$ is a closed analytic submanifold, the manifold $M_1$ is locally rigid in the sense that every nearby representation $G \to PU(2,1)$ stabilizes a complex geodesic in $\mathbb{H}_{\mathbb{C}}^2$ and is conjugate to a representation $G \to PU(1,1) \subset PU(2,1)$. In other words, there are no non-trivial "quasi-Fuchsian" deformations of the group $G$ and the complex surface $M_1$. On the other hand, as we show in Section 5 (cf. Theorem 5.1), the second manifold $M_2$ has plentiful enough Teichmüller space of different "quasi-Fuchsian" complex hyperbolic structures.

## 4. Deformations of holomorphic bundles: flexibility versus rigidity

Due to the natural inclusion $PO(n,1) \subset PU(n,1)$, any real hyperbolic $n$-manifold $M_{\mathbb{R}} = \mathbb{H}_{\mathbb{R}}^n/G$ can be (totally geodesically) embedded into a complex hyperbolic $n$-manifold $M_{\mathbb{C}} = \mathbb{H}_{\mathbb{C}}^n/G$ which is the total space of $n$-disk bundle over $M_{\mathbb{R}}$. Similarly, due to the inclusion $PU(n-1,1) \subset PU(n,1)$, any discrete torsion free group $G \subset PU(n-1,1)$ defines a holomorphic 2-disk bundle (with the total space $\mathbb{H}_{\mathbb{C}}^n/G$) over its totally geodesic complex analytic submanifolds $\mathbb{H}_{\mathbb{C}}^{n-1}/G$. In particular, one can consider both types of disk bundles over a Riemann surface $S$. Then a flexibility of such bundles becomes evident starting with hyperbolic structures on a Riemann surface $S$ of genus $g > 1$, which form the Teichmüller space $\mathcal{T}_g$, a complex analytic $(3g-3)$-manifold. And though, due to the Mostow rigidity theorem [Mo1], hyperbolic structures of finite volume and (real) dimension at least three are uniquely determined by their topology, so one has no continuous deformations of them, we still have some flexibility.

Firstly, real hyperbolic 3-manifolds have plentiful enough infinitesimal deformations and, according to Thurston's hyperbolic Dehn surgery theorem [T], noncompact hyperbolic 3-manifolds of finite volume can be approximated by compact hyperbolic 3-manifolds. Secondly, despite their hyperbolic rigidity, real hyperbolic manifolds $M$ can be deformed as conformal manifolds, or equivalently as higher-dimensional hyperbolic manifolds $M \times (0,1)$ of infinite volume. First such quasi-Fuchsian deformations were given by the author [A2] and, after Thurston's "Mickey Mouse" example [T], they were called bendings of $M$ along its totally geodesic hypersurfaces, see also [A1, A2, A4-A6, JM, Ko, Su1]. Furthermore, all these deformations are quasiconformally equivalent showing a rich supply of quasiconformal $G$-equivariant homeomorphisms in the real hyperbolic space $\mathbb{H}_{\mathbb{R}}^n$. In particular, the limit set $\Lambda(G) \subset \partial \mathbb{H}_{\mathbb{R}}^{n+1}$ deforms continuously from a round sphere $\partial \mathbb{H}_{\mathbb{R}}^n = S^{n-1} \subset S^n = \mathbb{H}_{\mathbb{R}}^{n+1}$ into nondifferentiably embedded topological $(n-1)$-spheres quasiconformally equivalent to $S^{n-1}$.

Contrasting to the above flexibility, "non-real" hyperbolic manifolds (locally symmetric spaces of rank one) seem much more rigid. In particular, due to P.Pansu [P], quasiconformal maps in the sphere at infinity of quaternionic/octanionic hyperbolic spaces that are induced by hyperbolic quasi-isometries are necessarily CR-automorphisms, and thus there cannot be interesting quasiconformal deformations of corresponding structures. Secondly, due to Corlette's rigidity theorem [Co2], such manifolds are even super-rigid – analogously to Margulis super-rigidity in higher rank [M]. The last fact implies impossibility of quasi-Fuchsian deformations of such quaternionic manifolds of dimension at least 3, see [Ka]. Furthermore, complex hyperbolic manifolds share the above rigidity of quaternionic/octanionic hyperbolic manifolds. Namely, due to the Goldman's local rigidity theorem in dimension $n = 2$ [G1] and its extension for $n \geq 3$ [GM], every nearby discrete representation $\rho : G \to PU(n, 1)$ of a cocompact lattice $G \subset PU(n-1, 1)$ stabilizes a complex totally geodesic subspace $\mathbb{H}_{\mathbb{C}}^{n-1}$ in $\mathbb{H}_{\mathbb{C}}^n$. Thus the limit set $\Lambda(\rho G) \subset \partial \mathbb{H}_{\mathbb{C}}^n$ is always a round sphere $S^{2n-3}$. Moreover, in higher dimensions $n \geq 3$, this local rigidity of complex hyperbolic $n$-manifolds $M$ homotopy equivalent to their closed complex totally geodesic hypersurfaces is even global due to a recent Yue's rigidity theorem [Y1]. These facts may be viewed as some arguments in favor of general rigidity and stability of deformations of complex hyperbolic structures.

To the contrary, our goal here and in the next section is to show that the opposite situation nevertheless holds: there are non-rigid complex hyperbolic manifolds which are disk bundles over their totally geodesic (both complex and real) submanifolds, and deformations of such manifolds may be quasiconformally unstable. The complex hyperbolic manifolds that are so flexible appear to be Stein spaces, so we expect that all Stein (complex hyperbolic) manifolds with "big" ends at infinity may have such nontrivial deformations.

The flexibility of complex hyperbolic 2-manifolds we deal with in this Section (and their property to be Stein spaces) is related to noncompactness of the (finite area) fibration base of these holomorphic disk bundles. In addition to constructing (quasi-Fuchsian) deformations of such bundles, we shall also show that they are quasiconformally unstable. Here we use the notion of quasiconformal stability that has its roots in the classical problem of quasiconformal stability of deformations from the theory of Kleinian groups, in particular in well known stability theorems by L. Bers [B1, B2] and D. Sullivan [Su1]. Let us recall that definition by following to L. Bers [B2]. Namely, a homomorphism $\chi : G \to \mathrm{PSL}(2, \mathbb{C})$ of a finitely generated group $G \subset \mathrm{PSL}(2, \mathbb{C})$ will be called allowable if it preserves the square traces of parabolic and elliptic elements (hence $\chi$ is type-preserving). A finitely generated Kleinian group $G \subset \mathrm{PSL}(2, \mathbb{C})$ is said to be *quasiconformally stable* if every allowable homomorphism $\chi : G \to \mathrm{PSL}(2, \mathbb{C})$ sufficiently close to the identity is induced by an equivarint quasiconformal mapping $w : \overline{\mathbb{C}} \to \overline{\mathbb{C}}$, that is $\chi(g) = wgw^{-1}$ for all $g \in G$. It is clear that degenerate Kleinian groups are quasiconformally unstable. However, due to a Bers's [B1] criterium (which involves quadratic differentials for the group $G$), it follows that Fuchsian groups, Schottky groups, groups of Schottky type and certain non-degenerate $B$-groups are all quasiconformally stable [B2].

We obtain a natural generalization of quasiconformal stability for discrete groups $G \subset \mathrm{Isom}\, \mathbb{H}_{\mathbb{C}}^n$ by changing the condition on homomorphisms $\chi$ in terms of the trace of elements $g \in G$ to the condition that such a homomorphism $\chi : G \to \mathrm{Isom}\, \mathbb{H}_{\mathbb{C}}^n$ preserves the type of elements of a given discrete group $G$. In that sense, B. Aebisher and R. Miner [AM] recently proved that (classical) Schottky groups $G \subset PU(n,1)$ are quasiconformally stable. Here a finitely generated discrete group $G = \langle g_1, \ldots, g_k \rangle \subset PU(n,1)$ is called a classical Schottky group of rank $k$ in the complex hyperbolic space $\mathbb{H}^n_\mathbb{C}$ if the sides of its Dirichlet polyhedron $D_y(G) \subset \mathbb{H}^n_\mathbb{C}$,

$$D_z(G) = \{z \in \mathbb{H}^n_\mathbb{C} : d(z,y) < d(z,g(y)) \text{ for any } g \in G\backslash\{\text{id}\}\}, \qquad (4.1)$$

centered at some point $y \in \mathbb{H}^n_\mathbb{C}$ are disjoint and non-asymptotic.

Nevertheless, as we shall show below, Fuchsian groups $G \subset PU(2,1)$ are quasiconformally unstable:

**Theorem 4.1.** *There are co-finite Fuchsian groups $G \subset PU(1,1) \subset PU(2,1)$ with signatures $(g, r; m_1, \ldots, m_r)$, where genus $g \geq 0$ and there are at least four cusps (with branching orders $m_i = \infty$), such that:*

(1) *the Teichmüller space $\mathcal{T}(G)$ contains a smooth simple curve $\alpha : [0, \pi/2) \hookrightarrow \mathcal{T}(g)$ which passes through the Fuchsian group $G = \alpha(0)$ and whose points $\alpha(t) = G_t \subset PU(2,1)$, $0 < t < \pi/2$, are all non-trivial quasi-Fuchsian groups;*
(2) *each isomorphism $\chi : G \to G_t$, $0 < t < \pi/2$, is induced by a $G$-equivariant homeomorphism $f_t : \overline{\mathbb{H}^2_\mathbb{C}} \to \overline{\mathbb{H}^2_\mathbb{C}}$ of the closure $\overline{\mathbb{H}^2_\mathbb{C}} = \mathbb{H}^2_\mathbb{C} \cup \partial \mathbb{H}^2_\mathbb{C}$ of the complex hyperbolic space;*
(3) *for any parameter $t$, $0 < t < \pi/2$, the action of the quasi-Fuchsian group $G_t$ is not quasiconformally conjugate to the action of the Fuchsian group $G = \alpha(0)$ (in both CR-structure at infinity $\partial \mathbb{H}^2_\mathbb{C} = \mathcal{H}_2 \cup \{\infty\}$ and the complex hyperbolic space $\mathbb{H}^2_\mathbb{C}$).*

Before we go on with the (constructive) proof of this Theorem, we note that though the constructed unstable Fuchsian groups $G \subset PU(1,1) \subset PU(2,1)$ of finite co-volume may have finite order elements, their finite index torsion free subgroups (and Riemann-Hurvitz formula for genus of a branching covering, see [KAG, (41)]) immediate imply the following:

**Corollary 4.2.** *Let $M = \mathbb{H}^2_\mathbb{C}/G$ be a complex hyperbolic surface with the holonomy group $G \subset PU(1,1) \subset PU(2,1)$ that represents the total space of a non-trivial disk bundle over a Riemann surface of genus $p \geq 0$ with at least four punctures (hyperbolic 2-orbifold with at least four punctures). Then the Teichmüller space $\mathcal{T}(M)$ contains a smooth simple curve $\alpha : [0, \pi/2) \hookrightarrow \mathcal{T}(M)$ with the following properties:*

(1) *the curve $\alpha$ passes through the surface $M = \alpha(0)$;*
(2) *each complex hyperbolic surface $M_t = \alpha(t) = \mathbb{H}^2_\mathbb{C}/G_t\}$, $t \in [0, \pi/2)$, with the holonomy group $G_t \subset PU(2,1)$ is homeomorphic to the surface $M$;*
(3) *for any parameter $t$, $0 < t < \pi/2$, the complex hyperbolic surface $M_t$ is not quasiconformally equivalent to the surface $M$.*

Besides the claims in this Corollary, it follows also from the construction of the above complex hyperbolic surfaces $M$ and $M_t$ that their boundaries, the spherical CR-manifolds $N = \partial M = \Omega(G)/G$ and $N_t = \partial M_t = \Omega(G_t)/G_t$ have similar properties:

**Corollary 4.3.** *Let $N = N_0 = \partial M$ be the 3-dimensional spherical CR-manifold with Fuchsian holonomy group $G \subset PU(1,1) \subset PU(2,1)$ that is the boundary at infinity of the complex hyperbolic surface $M$ from Corollary 4.2 (and which is the* *total space of a non-trivial circle bundle over a Riemann surface with at least four punctures). Then the Teichmüller space $\mathcal{T}(N)$ of the CR-manifold $N$ contains a smooth simple curve $\alpha_N : [0, \pi/2) \hookrightarrow \mathcal{T}(N)$ with the following properties:*

(1) *the curve $\alpha_N$ passes through the CR-manifold $N = \alpha_N(0)$;*
(2) *each CR-manifold $N_t = \alpha_N(t) = \mathbb{H}_{\mathbb{C}}^2/G_t\}$, $t \in [0, \pi/2)$, with the holonomy group $G_t \subset PU(2,1)$ is the total space of a non-trivial circle bundle over the Riemann surface with at least four punctures and is homeomorphic to the manifold $N$;*
(3) *for any parameter $t$, $0 < t < \pi/2$, the CR-manifold $N_t$ is not quasiconformally equivalent to the manifold $N$.*

*Remarks 4.4.*

(1) As a corollary of Theorem 4.1 and an Yue's [Y2] result on Hausdorff dimension, we have that there are deformations of a co-finite Fuchsian group $G \subset PU(1,1)$ into quasi-Fuchsian groups $G_\alpha = f_\alpha G f_\alpha^{-1} \subset PU(2,1)$ with Hausdorff dimension of the limit set $\Lambda(G_\alpha)$ strictly bigger than one. Moreover, the deformed groups $G_\alpha$ are Zariski dense in $PU(2,1)$.
(2) We note that, for the simplest case of manifolds with cyclic fundamental groups, a similar to Corollary 4.3 (though based on different ideas) effect of homeomorphic but not quasiconformally equivalent spherical CR-manifolds $N$ and $N'$ has been also recently observed by R. Miner [Mi]. In fact, among his Cauchy-Riemannian 3-manifolds (homeomorphic to $\mathbb{R}^2 \times S^1$), there are exactly two quasiconformal equivalence classes whose representatives have the cyclic holonomy groups generated correspondingly by a vertical Heisenberg translation by $(0,1) \in \mathbb{C} \times \mathbb{R}$ and a horizontal translation by $(1,0) \in \mathbb{C} \times \mathbb{R}$.
(3) The existence of non-trivial quasi-Fuchsian representations of Fuchsian groups with signatures $(g, r; m_1, \ldots, m_r)$ in the above Theorem 4.1 has its origin in an example (for genus $g = 0$ and $r = 4$ singular points with all four branching indices $m_i = \infty$) constructed by M. Carneiro and N. Gusevskii, see [Gu, CaG], who deal with a group $G \subset \text{Isom}\,\mathbb{H}_{\mathbb{C}}^1$ generated by four involutions and acting in the invariant complex geodesic in $\mathbb{H}_{\mathbb{C}}^2$ as the group generated by reflections in sides of an ideal 4-gon.
(4) Despite the impossibility of quasiconformal conjugation of the constructed actions of quasi-Fuchsian groups $G_t$ and $G = G_0$ in the sphere at infinity $\partial \mathbb{H}_{\mathbb{C}}^2$, it is still an open question whether the actions of these groups on their limit sets $\Lambda(G_t)$ and $\Lambda(G)$ could be "quasiconformally" conjugate, in other words, whether the canonical $G$-equivariant homeomorphism $f_{\chi_t}$, $f_{\chi_t} : \Lambda(G) \to \Lambda(G_t)$ of the limit sets (constructed in Theorem 3.1) that induces the isomorphism $\chi_t : G \to G_t$, $0 < t < \pi/2$, is in fact quasisymmetric.

*Proof of Theorem 4.1.* Here we present a complete construction of the deformation as well as basic arguments of the proof (for full details, see [A12]). We shall start with the lattice $G \subset PU(1,1)$ in the claim as a subgroup of index 2 in a discrete group $\Gamma \subset \text{Isom}\,\mathbb{H}_{\mathbb{C}}^1$ generated by reflections in sides of a finite area hyperbolic polygon $F \subset \mathbb{H}_{\mathbb{C}}^1 \subset \mathbb{H}_{\mathbb{C}}^2$.

Namely, let $\Gamma \subset \text{Isom}\,\mathbb{H}_{\mathbb{C}}^1$, $\Gamma = \Gamma_1 * \Gamma_2$, be a *co-finite* (free) lattice which is the free product of a dihedral parabolic subgroup $\Gamma_1$ and another subgroup $\Gamma_2$ which has at least one parabolic subgroup and whose quotient $[\overline{\mathbb{H}^1_{\mathbb{C}}}\backslash\Lambda(\Gamma_2)]/\Gamma_2$ has one boundary component at infinity, see Fig. 1. Since the subgroup $\Gamma_2$ has at least one parabolic subgroup, it can also be decomposed as the free product of its subgroups, $\Gamma_2 = \Gamma_3 * \Gamma_4$. In the simplest case, each of these subgroups $\Gamma_3$ and $\Gamma_4$ may have order two.

FIGURE 1. $\Gamma = \Gamma_1 * \Gamma_2$.

Here, for each reflection $g \subset \mathrm{Isom}\,\mathbb{H}^1_{\mathbb{C}}$ in a geodesic $\ell = (a,b) \subset \mathbb{H}^1_{\mathbb{C}}$, we have a uniquely defined (up to unitary rotation around the complex geodesic $\mathbb{H}^1_{\mathbb{C}} \subset \mathbb{H}^n_{\mathbb{C}}$) action of $g$ in the whole space $\mathbb{H}^n_{\mathbb{C}}$ as the anti-holomorphic involution whose fixed set is a real hyperbolic $n$-subspace in $\mathbb{H}^n_{\mathbb{C}}$ intersecting the complex geodesic $\mathbb{H}^1_{\mathbb{C}}$ along the geodesic $\ell = (a,b)$. In particular, assuming that $\mathbb{H}^1_{\mathbb{C}} = \{0\} \times \mathbb{R} \times \mathbb{R}_+ \subset \mathbb{C}^{n-1} \times \mathbb{R} \times \mathbb{R}_+ = \mathbb{H}^n_{\mathbb{C}}$ and that the geodesic $\ell$ ends at ideal points $a = (0,1), b = (0,-1) \in \mathcal{H}_n = \mathbb{C}^{n-1} \times \mathbb{R}$, we have that the reflection $g$ in $\ell$ acts in $\mathbb{H}^n_{\mathbb{C}}$ as the following anti-holomorphic involution (we call it a *real involution*):

$$i_\ell(\xi, v, u) = \left( \frac{A\overline{\xi}}{||\xi||^2 + u + iv}, \frac{v}{\left|||\xi||^2 + u + iv\right|^2}, \frac{u}{\left|||\xi||^2 + u + iv\right|^2} \right), \quad (4.2)$$

where $A \in U(n-1)$ and $(\xi, v, u) \in \mathbb{C}^{n-1} \times \mathbb{R} \times [0, \infty)$; compare with the Heisenberg inversion $\mathcal{I}$ in (2.1). The action of the real involution $i_\ell \in \mathrm{Isom}\,\mathbb{H}^n_{\mathbb{C}}$ at infinity $\partial\mathbb{H}^n_{\mathbb{C}} = \mathbb{C}^{n-1} \times \mathbb{R} \times \{0\} \cup \{\infty\}$ preserves the unit Heisenberg sphere $S_c(0,1) = \{(\xi,v) \in \mathcal{H}_n : ||(\xi,v)||_c = 1\}$, swaps the origin $(0,0) \in \mathcal{H}_n$ and $\infty$, and pointwise fixes an $\mathbb{R}$-sphere (of dimension $(n-1)$) that lies in the Heisenberg sphere $S_c(0,1)$ and passes through its poli $a = (0,1), b = (0,-1) \in \mathcal{H}_n = \mathbb{C}^{n-1} \times \mathbb{R}$. In particular, for $n=2$ and $A\xi = -\xi$, the pointwise fixed $\mathbb{R}$-circle has the following equation in cylindrical coordinates (cf. [G4]):

$$\{(re^{i\theta}, v) \in \mathcal{H} = \mathbb{C} \times \mathbb{R} : r^2 + iv = -e^{2i\theta}, \quad (4.3)$$

and its vertical projection to the horizontal plane $\mathbb{C} \times \{0\} \subset \mathcal{H}$ is the lemniscate of Bernoulli, see Fig.2:

$$\{\xi = x + iy : (x^2 + y^2)^2 + x^2 - y^2 = 0\}.$$

FIGURE 2. Vertical projection of $\mathbb{R}$-circle.

Now we may assume that our group $\Gamma \subset \operatorname{Isom} \mathbb{H}_{\mathbb{C}}^1$ is generated by real involutions whose restrictions to $\mathbb{H}_{\mathbb{C}}^1$ are reflections in sides of the fundamental polygon $F$, and its limit set is the (vertical) chain $\{0\} \times \mathbb{R} \cup \{\infty\} \subset \partial \mathbb{H}_{\mathbb{C}}^2$. Furthermore, assuming that the (parabolic) fixed point of the dihedral subgroup $\Gamma_1 \subset \Gamma$ is $\infty$ and deforming the action of $\Gamma$ in $\mathbb{H}_{\mathbb{C}}^1$ (in Teichmüller space $\mathcal{T}(\Gamma)$), we may take a fundamental polyhedron $D \subset \mathbb{H}_{\mathbb{C}}^2$ of the group $\Gamma$ as the polyhedron bounded by bisectors $\Sigma_i$ whose poli lie in the one point compactification of the vertical line $\{0\} \times \mathbb{R} \subset \mathcal{H}$ and whose boundaries $S_i$ at infinity are as follows. Two of them (corresponding to the generators of $\Gamma_1$) are the extended horizontal planes in $\mathcal{H}$,

$$S_1 = \mathbb{C} \times \{s_0\} \cup \{\infty\}, \ S_2 = \mathbb{C} \times \{-s_0\} \cup \{\infty\} \subset \mathbb{C} \times \mathbb{R} \cup \{\infty\}. \qquad (4.4)$$

The spheres at infinity of all other bisectors are Heisenberg spheres $S_i$, $i \geq 3$, that lie between the planes (4.4) and whose centers lie in the vertical line. Furthermore, we may assume that two such spheres, $S_3$ and $S_4$, are unit Heisenberg spheres tangent to the spheres $S_1$ and $S_2$, correspondingly at the points $p_1 = (0, s_0)$ and $p_2 = (0, -s_0)$. Also, due to our condition on one more parabolic subgroup (in $\Gamma_2$), we have that there is a pair of spheres, $S_j$ and $S_{j+1}$, $j \geq 3$, tangent to each other at some point $p_3 = (0, s_1)$. Figure 3 shows such a configuration of (at least four) spheres $\{S_i\}$ as well as a choice of $\mathbb{R}$-circles on them, $m_i \subset S_i$, that are fixed by the corresponding real involutions $\gamma_i \in \operatorname{Isom} \mathbb{H}_{\mathbb{C}}^2$ that generate the group $\Gamma$. We note that those involutions $\gamma_i$ are lifts of reflections in sides of the polygon $F$ in the complex geodesic $\mathbb{H}_{\mathbb{C}}^1 \subset \mathbb{H}_{\mathbb{C}}^2$. So those lifts should be compatible in the sense that the union of closures of real arcs $m_i \cap \overline{D}$ in their pointwise fixed $\mathbb{R}$-circles $m_i$ is a closed loop $\alpha$ on the boundary of the 3-dimensional polyhedron $P = \overline{D} \cap \partial \mathbb{H}_{\mathbb{C}}^2$.

FIGURE 3. Fundamental polyhedron $P$ and generators of $\Gamma$.

As such $\mathbb{R}$-circles $m_i \subset S_i$, $i = 1, 2$, we take two real (extended) lines that are parallel to each other and pass the corresponding poli $p_i = (0, \pm s_0)$ in $S_i$. In the adjacent spheres $S_j$, $j = 3, 4$, we take $\mathbb{R}$-circles $m_j$ as those unique circles that are tangent at the points $(0, \pm s_0)$ to to corresponding $\mathbb{R}$-circles $m_1$ and $m_2$. Continuing this process in the spheres $S_l$ adjacent to $S_3$ and $S_4$, we take those $\mathbb{R}$-circles $m_l \subset S_l$ that intersect the corresponding $\mathbb{R}$-circles $m_i \subset S_i$, $i = 3, 4$. Continuing this process, we will reach next tangent points, etc. Finally, at a tangent point $(0, s_1)$, our $\mathbb{R}$-circles $m_j$ and $m_j + 1$ will meet at some angle $\theta_{j,j+1}$. Then our group $\Gamma$,

$$\Gamma = \langle \gamma_1, \gamma_2, \gamma_3, \ldots, \gamma_k \rangle = \langle \gamma_1 \rangle * \langle \gamma_2 \rangle * \Gamma_3 * \Gamma_4 \,, \qquad (4.5)$$

is generated by the real involutions $\gamma_i$ pointwise fixing the corresponding $\mathbb{R}$-circles $m_i$, and acts in the invariant complex geodesic $\mathbb{H}_{\mathbb{C}}^1 \subset \mathbb{H}_{\mathbb{C}}^2$ as the group generated by reflections in sides of hyperbolic $k$-gon $P_0$ of finite area. Also we note that in our choice of lifts of generating reflections for each of the subgroups $\Gamma_3, \Gamma_4 \subset \Gamma_2$ as real involutions in $\mathbb{H}_{\mathbb{C}}^2$, we have some degree of freedom (compositions of those involutions with unitary rotations about the vertical line $\{0\} \times \mathbb{R} \subset \mathcal{H}$, see (4.2)) which depends on the number of conjugacy classes of parabolic subgroups in $\Gamma_2$.

Obviously, we have that the subgroup $G \subset \Gamma$ of index two with the fundamental polyhedron $\overline{D} \cup \gamma_1(\overline{D})$ is a subgroup of $PU(1,1) \subset PU(2,1)$. Topologically, the quotient $\mathbb{H}_{\mathbb{C}}^2/G$ is a 2-disk bundle over 2-dimensional sphere with at least four punctures; geometrically, the base of this bundle is its totally geodesic complex analytic orbifold of finite hyperbolic volume (due to Riemann-Hurvitz formula for genus of a branching covering over 2-sphere, see [KAG, (41)], it is covered by a Riemann surface of a genus $p \geq 0$ with at least four punctures).

Also we note that tangent points of bisectors bounding the fundamental polyhedron $D \subset \mathbb{H}_\mathbb{C}^2$, in particular the points $p_0 = \infty$ and $p_i = (0, \pm s_0) \in \mathcal{H}$, $i = 1, 2$, are parabolic fixed points of $\Gamma$. Moreover, due to tangency of the corresponding $\mathbb{R}$-circles in the definition of $\Gamma$, all elements $g \in G_{p_i}$ in the stabilizer subgroups $G_{p_i} \subset G \subset PU(2,1)$ of those last three points are (conjugated to) pure vertical translations, that is (after conjugation) they act in $\overline{\mathbb{H}_\mathbb{C}^2}$ as Heisenberg translations $(\xi, v, u) \mapsto (\xi, v + v_i, u)$, $i = 0, 1, 2$.

**Deformation of groups $\Gamma$ and $G \subset \Gamma$.**

To deform the groups $\Gamma$ and $G \subset \Gamma$, we define a family of discrete faithfull representations $\rho_t : \Gamma \to \mathrm{Isom}\,\mathbb{H}_\mathbb{C}^2$, $0 \leq t < \pi/2$, with the images $\rho(t) = \Gamma^t$ where $\rho(0) = \Gamma^0 = \Gamma$.

Namely, all these representations $\rho_t$ coincide (up to conjugations by unitary rotations $A$ in (4.2)) on the subgroup $\Gamma_2 \subset \Gamma$ and only differ on the dihedral subgroup $\Gamma_1 \subset \Gamma$ in the following way:

$$\Gamma^t = \langle \gamma_{1,t}, \gamma_{2,t}, \gamma_{3,t} \ldots, \gamma_{k,t} \rangle = \langle \gamma_{1,t} \rangle * \langle \gamma_{2,t} \rangle * \Gamma_3 * A_t \Gamma_4 A_t^{-1}, \qquad (4.6)$$

where $A_t \in U(1)$ acts in $\mathbb{H}_\mathbb{C}^2$ by unitary rotation about the complex geodesic $\mathbb{H}_\mathbb{C}^1 \subset \mathbb{H}_\mathbb{C}^2$, and the generators $\gamma_{i,t}$ are anti-holomorphic (real) involutions with pointwise fixed $\mathbb{R}$-circles $m_{i,t} \subset \partial \mathbb{H}_\mathbb{C}^2$. In particular, $\gamma_{1,t}$ and $\gamma_{2,t}$ generate the new dihedral subgroup $\Gamma_{1,t} = \rho_t(\Gamma_1)$ of the group $\Gamma^t \subset \mathrm{Isom}\,\mathbb{H}_\mathbb{C}^2$.

As before, we define the real involutions $\gamma_{i,t}$ by determining their fixed $\mathbb{R}$-circles $m_{i,t} \subset \partial \mathbb{H}_\mathbb{C}^2$. Namely, for each $t$, $0 \leq t < \pi/2$, let $p_{1,t}$ be a point on the $\mathbb{R}$-circle $m_3$ that is seen from the center of the sphere $S_3$ at the angle $t$, and let $p_{2,t}$ be the point on the Heisenberg sphere $S_4$ that is symmetric to the point $p_{1,t}$, $p_{2,t} = -p_{1,t}$. Then let $m_{4,t}$ be the $\mathbb{R}$-circle in the sphere $S_4$ that contains the point $p_{2,t}$, and $A_t \in U(1)$ be the unitary rotation of $S_4$ that maps $m_4$ to $m_{4,t}$. We note that $A_t \neq \mathrm{id}$ if $t \neq 0$, and its rotation angle monotonically increases from 0 to $\pi/2$ as $t$ tends from 0 to $\pi/2$. Now we keep all generators in the subgroup $\Gamma_3$ and conjugate generators of the group $\Gamma_4$ by $A_t$. The remaining $\mathbb{R}$-circles $m_{1,t}, m_{2,t} \subset \partial \mathbb{H}_\mathbb{C}^2$ are given as the (Euclidean) lines in $\mathbb{C} \times \mathbb{R}$ tangent to the corresponding $\mathbb{R}$-circles $m_{3,t} = m_3$ and $m_{4,t} = A_t(m_4)$ at the points $p_{1,t}$ and $p_{2,t}$, respectively.

In other words, we replace the bisectors $\Sigma_1$ and $\Sigma_2$ by bisectors $\Sigma_1^t$ and $\Sigma_2^t$ whose boundaries at infinity are the (extended) parallel planes $S_1^t$ and $S_2^t$ tangent to the spheres $S_3$ and $S_4$ at the points $p_{1,t}$ and $p_{2,t}$, respectively. Then the previously defined $\mathbb{R}$-circles $m_{1,t}$ and $m_{2,t}$ are the intersection lines of the planes $S_1^t$ and $S_2^t$ with the contact planes at these tangent points $p_{1,t}$ and $p_{2,t}$, correspondingly, see Fig. 4. Clearly, these lines are not parallel if $t \neq 0$, and the angle between them (at infinity) monotonically increases from 0 to $\pi/2$ as $t$ tends from 0 to $\pi/2$. It is also worth to mention that, for $t \neq 0$, the finite poli of the spheres $S_1^t$ and $S_2^t$ lie somewhere in the $\mathbb{R}$-circles $m_{1,t}$ and $m_{2,t}$ and differ from the tangent points $p_{1,t}$ and $p_{2,t}$, respectively.

We note that as in the group $\Gamma$, all tangent points of bisectors $\Sigma_i^t$ and $\Sigma_j^t$ (and in particular, the points $p_1^t, p_2^t$ and $\infty$) are parabolic fixed points of the deformed group $\Gamma^t$. Moreover, though parabolic elements fixing the points $p_1^t$ and $p_2^t$ are still conjugate to vertical Heisenberg translations, the index two subgroup of projective unitary transformations in the stabilizer of $\infty$ in $\Gamma^t$ is generated by screw vertical translation if $t \neq 0$. The rotation angle of that screw translation monotonically increases from 0 to $\pi$ as $t$ tends from 0 to $\pi/2$.

FIGURE 4. Fundamental polyhedron $P^t$ and generators of $\Gamma^t$.

The above property of tangent points implies that the polyhedron $D^t \subset \mathbb{H}^2_{\mathbb{C}}$ bounded by bisectors $\Sigma^t_i$ is a fundamental polyhedron of the discrete group $\Gamma^t$, and its intersection with the sphere at infinity, the polyhedron $P^t = \overline{D^t} \cap \partial\mathbb{H}^2_{\mathbb{C}}$ bounded by the spheres $S^t_i$, is a fundamental polyhedron for the $\Gamma^t$-action at infinity shown in Fig. 4.

As another implication of properties of parabolic subgroups in $\Gamma$ and $\Gamma^t$, we obtain non-triviality of our deformation given by the family of (faithful discrete) representations $\rho_t \in \mathrm{Hom}(\Gamma, \mathrm{Isom}\,\mathbb{H}^2_{\mathbb{C}})$, $\rho_t \Gamma = \Gamma^t$. Namely, for any two different parameters $t$ and $t'$, the groups $\Gamma^t$ and $\Gamma^{t'}$ cannot be conjugated in $\mathrm{Isom}\,\mathbb{H}^2_{\mathbb{C}}$ because the corresponding parabolic transformations in their stabilizers of $\infty$,

$$\rho_t(\gamma) \in \Gamma^t_\infty \text{ and } \rho_{t'}(\gamma) \in \Gamma^{t'}_\infty, \tag{4.7}$$

have rotational parts (unitary rotations) with different angles. We note that one may also derive this fact by using the Cartan angular invariant for triples of points $((0, -s_0), (0, s_1), (0, s_0))$ and $(p_{2,t}, (0, s_1), p_{1,t})$. This Cartan invariant (see the next Section) is different from $\pm\pi/2$ for any $t \neq 0$, which also shows that the groups $\Gamma^t$ are non-trivial quasi-Fuchsian groups (with the limit topological circles $\Lambda(\Gamma^t)$ different from "round" circles in $\partial\mathbb{H}^2_{\mathbb{C}}$). Nevertheless, these quasi-Fuchsian groups $\Gamma^t$ are not quasiconformal conjugates of the Fuchsian group $\Gamma$:

**Lemma 4.5.** *For any parameter $t$, $0 < t < \pi/2$, the quasi-Fuchsian representation $\rho_t$ cannot be conjugate to the (Fuchsian) inclusion $\rho_0 : \Gamma \subset PU(2,1)$ by any quasiconformal homeomorphism in neither $\mathbb{H}^2_{\mathbb{C}}$ nor the Heisenberg group $\mathcal{H}$.*

*Proof.* We can prove this fact by using the above observation that the (type preserving) isomorphism $\rho_t : \Gamma \to \Gamma^t$ strictly increases the angles of rotational parts of parabolic elements in (4.7).

Namely, let us consider a family $F$ of vertical intervals in $\mathcal{H}$ connecting points $x$ in the plane $S_2$ to their images $\gamma(x)$ in the parallel plane $\gamma(S_2)$, with respect to the (pure) Heisenberg translation $\gamma = \gamma_1 \gamma_2$. Then, for any $\Gamma$-equivariant homeomorphism $f_t$ in the Heisenberg group $\mathcal{H}$ that conjugates $\Gamma$ to $\Gamma^t$, this family $F$ is mapped to a family $f_t(F)$ consisting of curves that connect the points $f_t(x)$ to their images $\gamma_t(f_t(x))$ with respect to the (screw) Heisenberg translation $\gamma_t = \gamma_{1,t} \gamma_{2,t}$ whose angle of rotation is not zero.

Due to Theorem 2.1, the stabilizer of $\infty$ in the deformed group $\Gamma^t$ has an invariant vertical line $\ell \subset \mathcal{H}$. So, for any point $f_t(x)$ in the plane $S_2^t$ having the distance $r$ to the line $\ell$, its image $\gamma_t(f_t(x))$ in the plane $\gamma_t(S_2^t)$ has the same distance $r$ to the line $\ell$. However when the distance $r$ becomes arbitrarily large, the length of the curve in the family $f_t(F)$ that connects the points $f_t(x)$ and $\gamma_t(f_t(x))$ goes to $\infty$, despite the fact that the lenght of the verical interval connecting the points $x$ and $\gamma(x)$ is constant. Due to the fact that any $K$-quasiconformal map maps curves of bounded length to curves whose lengths are bounded as well, the above observation shows that any such a $\Gamma$-equivariant homeomorphism $f_t$ in the Heisenberg group $\mathcal{H}$ cannot be quasiconformal. This completes the proof of Lemma 4.5.

$\square$

We shall finish the proof of Theorem 4.1 by showing that our deformaton is nevertheless topological, that is it could be induced by a continuous family of equivariant homeomorphisms:

**Proposition 4.6.** *For any two parameters $t$ and $t'$, $0 \leq t < t', \pi/2$, the faithful discrete (quasi-Fuchsian) representations $\rho_t : \Gamma \to \Gamma^t$ and $\rho_{t'} : \Gamma \to \Gamma^{t'}$ are conjugate by an equivariant homeomorphism (that continuously depends on $t$ and $t'$),*

$$f_{t,t'} : \overline{\mathbb{H}^2_{\mathbb{C}}} \to \overline{\mathbb{H}^2_{\mathbb{C}}}. \tag{4.8}$$

*Proof.* It is enough to construct such an equivariant homeomorphism $f_t = f_{0,t}$ in (4.8) that conjugates the groups $\Gamma$ and $\Gamma^t$ for any $t \in [0, \pi/2)$ and continuously depends on $t$. To do that, we need to define such an equivariant homeomorphism $f_t$ on the closure $\overline{D} \subset \overline{\mathbb{H}^2_{\mathbb{C}}}$ of the fundamental polyhedron $D \subset \mathbb{H}^2_{\mathbb{C}}$ of the group $\Gamma$, whose image $f_t(\overline{D}) = \overline{D^t}$ is the closure of the fundamental polyhedron of the deformed group $\Gamma^t$. Having such a map of closed polyhedra, we can immediately extend it equivariantly to the whole discontinuity domain, that is to a $\Gamma$-equivariant homeomorphism

$$f_t : \mathbb{H}^2_{\mathbb{C}} \cup \Omega(\Gamma) \to \mathbb{H}^2_{\mathbb{C}} \cup \Omega(\Gamma^t), \tag{4.9}$$

whose extension by continuity to the limit set is the unique $\Gamma$-equivariant homeomorphism of the limit sets, $\Lambda(\Gamma) \to \Lambda(\Gamma^t)$, induced by the type preserving isomorphism of the groups $\Gamma$ and $\Gamma^t$, see Theorem 3.2.

Now we start a construction of our $\Gamma$-equivariant homeomorphism $f_t : \overline{D} \to \overline{D^t}$ on the closure of the fundamental polyhedron $D$, where it maps $k$-sides of $D$ to the corresponding $k$-sides of $D^t$, $0 \leq k \leq 3$, in particular tangent points $p_i$ to tangent points $p_{i,t}$, and the $\mathbb{R}$-circles $m_j$ to the corresponding $\mathbb{R}$-circles $m_{j,t}$.

First, we define our homeomorphism $f_t$ in the 3-polyhedron $P \subset \mathcal{H} \subset \partial \mathbb{H}_{\mathbb{C}}^2$ at infinity. Fixing an orientation on the $\mathbb{R}$-circles $m_j$ and $m_{j,t}$ (preserved by their identifications), we consider corresponding (positive) semispheres $\hat{S}_j \subset S_j$ and $\hat{S}_j^t \subset S_j^t$ bounded by those $\mathbb{R}$-circles. Since, up to $PU(2,1)$, each such a semisphere may be considered as a halfplane in an extended 2-plane, we define homeomorphisms $f_t|_{\hat{S}_j}$ as restrictions of either Euclidean isometries in $\mathbb{R}^3 = \mathbb{C} \times \mathbb{R}$ that map $\hat{S}_j \to \hat{S}_j^t$ or conjugations of such isometries by elements of $PU(2,1)$. On the complements $S_j \backslash \hat{S}_j$, we define $f_t$ so that they are compatible with the generators of the group $\Gamma$:

$$f_t|_{S_j \backslash \hat{S}_j}(x) = \gamma_{j,t} \circ f_t|_{S_j} \circ \gamma_j(x), \quad \text{for } x \in S_j \backslash \hat{S}_j. \tag{4.10}$$

It defines our $\Gamma$-equivariant homeomorphism $f_t$ on the boundary $\partial P$ of the fundamental polyhedron $P \subset \mathcal{H}$, $P = \partial_\infty D$, which is not simply-connected.

Now we consider a closed curve $\alpha \subset \partial P$ as the union of arcs of our $\mathbb{R}$-circles $m_j$. In particular, this curve $\alpha$ connects $\infty$ and all other tangent (parabolic fixed) points in $\partial P$, and bounds an embedded topological 2-disk $\sigma \subset P$, $\partial \sigma = \alpha$. The disk $\sigma$ splits the polyhedron $P$ into a topological 3-ball $\hat{P}$, see Fig.5.

FIGURE 5. Cell decomposition of fundamental polyhedron $P$.

Similarly, on the boundary of the deformed polyhedron $P^t = \partial_\infty D^t$, we consider a closed curve $\alpha_t \subset \partial P^t$ as the union of arcs of our $\mathbb{R}$-circles $m_{j,t} = f_t(m_j)$, that is $\alpha_t = f_t(\alpha)$. As before, let $\sigma_t \subset \text{int } P_t$ be an embedded topological 2-disk spanned by the loop $\alpha_t$, that splits the polyhedron $P$ into a topological 3-ball $\hat{P}^t$. Then we can define our homeomorphism $f_t$ on the cell $\sigma$ as a homeomorphism $\sigma \to \sigma_t$ extending the homeomorphism $\alpha \to \alpha_t$ between the boundary loops.

Now we can cut the polyhedra $P$ and $P^t$ along the disks $\sigma$ and $\sigma_t$ correspondingly, and obtain two topological 3-balls $\hat{P}$ and $\hat{P}^t$ whose boundary spheres are $\Gamma$-equivariantly homeomorphic under our constructed map $f_t$. Extending that boundary map to a homeomorphism of the closed 3-balls, we finally have defined our $\Gamma$-equivariant homeomorphism $f_t|_P : P \to P^t$.

In addition to 3-polyhedra $P$ and $P^t$, the boundaries of the fundamental 4-dimensional polyhedra $D$ and $D^t$ each have $2k$ more sides. Those sides lie on bisectors $\Sigma_j$ and $\Sigma_j^t$ and are pairwise identified by the generators $\gamma_j$ and $\gamma_{j,t}$, respectively. Each of those bisectors is split into two halves, $\Sigma_j^+$ and $\Sigma_j^-$, or $\Sigma_j^{t,+}$ and $\Sigma_j^{t,-}$, along the corresponding totally real geodesic 2-plane spanned by the corresponding $\mathbb{R}$-circle, either $m_j$ or $m_{j,t}$, and those halves are pairwise identified by our generators.

Since we have already defined our homeomorphism $f_t$ on the boundary of each bisector $\Sigma_j$, we can extend it to the whole bisector by using natural foliations of $\Sigma_j$ and $\Sigma_j^t$ by disjoint real geodesics with ends $x$, $\gamma_j(x)$ and $f_t(x)$, $\gamma_{j,t}f_t(x)$, respectively. Namely, on each half $\overline{\Sigma_j^+}$, we define our homeomorphism $f_t$ as the map $\overline{\Sigma_j^+} \to \overline{\Sigma_j^{t,+}}$ that maps each (real) geodesic ray $(x, \gamma_j(x)) \cap \Sigma_j^+$ to the geodesic ray $(f_t(x), \gamma_{j,t}f_t(x)) \cap \Sigma_j^{t,+}$, $j = 1, 2, \ldots, k$. In the other halves, $\Sigma_j^-$, we define $f_t$ compatible with the generators of our groups (cf. (4.10)):

$$f_t|_{\Sigma_j^-}(x) = \gamma_{j,t} \circ f_t|_{\Sigma_j^+} \circ \gamma_j(x), \quad \text{for } x \in \Sigma_j^+. \tag{4.11}$$

Now, having our $\Gamma$-equivariant homeomorphism $f_t$ defined on the topological boundary 3-sphere $\partial D$,

$$f_t|_{\partial D} : \partial D \to \partial D^t,$$

we can extend this homeomorphism to a homeomorphism of the closed 4-balls, $\overline{D} \to \overline{D^t}$, that conjugates the dynamics of the groups $\Gamma$ and $\Gamma^t$ on the closures of our fundamental polyhedra. Clearly, all the steps in our construction of those homeomorphisms continuously depend on $t$. Hence, the equivariant extension (4.9) of that map to the whole 4-space completes the proof of our Proposition as well as of the whole Theorem 4.1. □

## 5. Bending deformations of complex hyperbolic structures

Here we present another class of non-rigid complex hyperbolic manifolds fibering over a Riemann surface. In fact we point out that the non-compactness condition for the base of fibration in non-rigidity results in §4 is not essential, either. Namely, complex hyperbolic Stein manifolds $M$ homotopy equivalent to their closed totally *real* geodesic surfaces are not rigid, too. To prove that, we shall present a canonical construction of continuous non-trivial quasi-Fuchsian deformations of complex surfaces fibered over closed Riemannian surfaces of genus $g > 1$. Such deformations depend on $3(g-1)$ real-analytic parameters (in addition to "Fuchsian" deformations, where in particular, the Teichmuller space of the base surface has dimension $6(g-1)$). This provides the first such non-trivial deformations of fibrations with compact base (for non-compact base, see a different Goldman-Parker' deformation [GP2] of ideal triangle groups $G \subset PO(2,1)$). The obtained flexibility of such holomorphic fibrations and the number of its parameters provide a partial confirmation of a conjecture on dimension $16(g-1)$ of the Teichmuller space of such complex surfaces. It is related to A.Weil's theorem [W] (see also [G3, p.43]), that the variety of conjugacy classes of all (not necessarily discrete) representations $G \to PU(2,1)$ near the embedding $G \subset PO(2,1)$ is a real-analytic manifold of dimension $16(g-1)$. We remark that discreteness of representations of $G \cong \pi_1 M$ is an essential condition for deformation of a complex manifold $M$ which does not follow from the mentioned Weil's result.

Our construction is inspired by the approach the author used for bending deformations of real hyperbolic (conformal) manifolds along totally geodesic hypersurfaces ([A2, A4]). In the case of complex hyperbolic (and Cauchy-Riemannian) structures, the constructed "bendings" work however in a different way than in the real hyperbolic case. Namely our complex bending deformations involve simultaneous bending of the base of the fibration of the complex surface $M$ as well as bendings of each of its totally geodesic fibers (see Remark 5.9). Such bending deformations of complex surfaces are associated to their real simple closed geodesics (of real codimension 3), but have nothing common with the so called cone deformations of real hyperbolic 3-manifolds along closed geodesics, see [A6, A9].

Furthermore, despite well known complications in constructing equivariant homeomorphisms in the complex hyperbolic space and in Cauchy-Riemannian geometry (which should preserve Kähler and contact structures in $\mathbb{H}_{\mathbb{C}}^n$ and at its infinity $\overline{\mathcal{H}_n}$, respectively), the constructed complex bending deformations are induced by equivariant homeomorphisms of $\overline{\mathbb{H}_{\mathbb{C}}^n}$. Moreover, in contrast to the situation with deformations in §4, those equivariant homeomorphisms are in addition quasiconformal:

**Theorem 5.1.** *Let $G \subset PO(2,1) \subset PU(2,1)$ be a given (non-elementary) discrete group. Then, for any simple closed geodesic $\alpha$ in the Riemann 2-surface $S = H_{\mathbb{R}}^2/G$ and a sufficiently small $\eta_0 > 0$, there is a holomorphic family of $G$-equivariant quasiconformal homeomorphisms $F_\eta : \overline{\mathbb{H}_{\mathbb{C}}^2} \to \overline{\mathbb{H}_{\mathbb{C}}^2}$, $-\eta_0 < \eta < \eta_0$, which defines the bending (quasi-Fuchsian) deformation $\mathcal{B}_\alpha : (-\eta_0, \eta_0) \to \mathcal{R}_0(G)$ of the group $G$ along the geodesic $\alpha$, with $\mathcal{B}_\alpha(\eta) = F_\eta^*$.*

We notice that deformations of a complex hyperbolic manifold $M$ may depend on many parameters described by the Teichmüller space $\mathcal{T}(M)$ of isotopy classes of complex hyperbolic structures on $M$. One can reduce the study of this space $\mathcal{T}(M)$ to studying the variety $\mathcal{T}(G)$ of conjugacy classes of discrete faithful representations $\rho : G \to PU(n,1)$ (involving the space $\mathcal{D}(M)$ of the developing maps, see [G2, FG]). Here $\mathcal{T}(G) = \mathcal{R}_0(G)/PU(n,1)$, and the variety $\mathcal{R}_0(G) \subset \text{Hom}(G, PU(n,1))$ consists of discrete faithful representations $\rho$ of the group $G$ with infinite co-volume, $\text{Vol}(\mathbb{H}_{\mathbb{C}}^n/G) = \infty$. In particular, our complex bending deformations depend on many independent parameters as it can be shown by applying our construction and Élie Cartan [Car] angular invariant in Cauchy-Riemannian geometry:

**Corollary 5.2.** *Let $S_p = \mathbb{H}_{\mathbb{R}}^2/G$ be a closed totally real geodesic surface of genus $p > 1$ in a given complex hyperbolic surface $M = \mathbb{H}_{\mathbb{C}}^2/G$, $G \subset PO(2,1) \subset PU(2,1)$. Then there is a real-analytic embedding $\pi \circ \mathcal{B} : B^{3p-3} \hookrightarrow \mathcal{T}(M)$ of a real $(3p-3)$-ball into the Teichmüller space of $M$, defined by bending deformations along disjoint closed geodesics in $M$ and by the projection $\pi : \mathcal{D}(M) \to \mathcal{T}(M) = \mathcal{D}(M)/PU(2,1)$ in the development space $\mathcal{D}(M)$.*

**Basic Construction (Proof of Theorem 5.1).** Now we start with a totally

real geodesic surface $S = \mathbb{H}^2_\mathbb{R}/G$ in the complex surface $M = \mathbb{H}^2_\mathbb{C}/G$, where $G \subset PO(2,1) \subset PU(2,1)$ is a given discrete group, and fix a simple closed geodesic $\alpha$ on $S$. We may assume that the loop $\alpha$ is covered by a geodesic $A \subset \mathbb{H}^2_\mathbb{R} \subset \mathbb{H}^2_\mathbb{C}$ whose ends at infinity are $\infty$ and the origin of the Heisenberg group $\mathcal{H} = \mathbb{C} \times \mathbb{R}$, $\overline{\mathcal{H}} = \partial \mathbb{H}^2_\mathbb{C}$. Furthermore, using quasiconformal deformations of the Riemann surface $S$ (in the Teichmüller space $\mathcal{T}(S)$, that is, by deforming the inclusion $G \subset PO(2,1)$ in $PO(2,1)$ by bendings along the loop $\alpha$, see Corollary 3.3 in [A10]), we can assume that the hyperbolic length of $\alpha$ is sufficiently small and the radius of its tubular neighborhood is big enough:

**Lemma 5.3.** *Let $g_\alpha$ be a hyperbolic element of a non-elementary discrete group $G \subset PO(2,1) \subset PU(2,1)$ with translation length $\ell$ along its axis $A \subset \mathbb{H}^2_\mathbb{R}$. Then any tubular neighborhood $U_\delta(A)$ of the axis $A$ of radius $\delta > 0$ is precisely invariant with respect to its stabilizer $G_0 \subset G$ if $\sinh(\ell/4) \cdot \sinh(\delta/2) \leq 1/2$. Furthermore, for sufficiently small $\ell$, $\ell < 4\delta$, the Dirichlet polyhedron $D_z(G) \subset \mathbb{H}^2_\mathbb{C}$ of the group $G$ centered at a point $z \in A$ has two sides $a$ and $a'$ intersecting the axis $A$ and such that $g_\alpha(a) = a'$.*

Then the group $G$ and its subgroups $G_0, G_1, G_2$ in the free amalgamated (or HNN-extension) decomposition of $G$ have Dirichlet polyhedra $D_z(G_i) \subset \mathbb{H}^2_\mathbb{C}$, $i = 0, 1, 2$, centered at a point $z \in A = (0, \infty)$, whose intersections with the hyperbolic 2-plane $\mathbb{H}^2_\mathbb{R}$ have the shapes indicated in Figures 6-9.

FIGURE 6.  $G_1 \subset G = G_1 *_{G_0} G_2$     FIGURE 7.  $G_2 \subset G = G_1 *_{G_0} G_2$

FIGURE 8.  $G_1 \subset G = G_1 *_{G_0}$     FIGURE 9.  $G = G_1 *_{G_0}$

In particular we have that, except of two bisectors $\mathfrak{S}$ and $\mathfrak{S}'$ that are equivalent under the hyperbolic translation $g_\alpha$ (that generates the stabilizer $G_0 \subset G$ of the axis $A$), all other bisectors bounding such a Dirichlet polyhedron lie in sufficiently small "cone neighborhoods" $C_+$ and $C_-$ of the arcs (infinite rays) $\mathbb{R}_+$ and $\mathbb{R}_-$ of the real circle $\mathbb{R} \times \{0\} \subset \mathbb{C} \times \mathbb{R} = \mathcal{H}$.

Actually, we may assume that the Heisenberg spheres at infinity of the bisectors $\mathfrak{S}$ and $\mathfrak{S}'$ have radii 1 and $r_0 > 1$, correspondingly. Then, for a sufficiently small $\epsilon$, $0 < \epsilon << r_0 - 1$, the cone neighborhoods $C_+, C_- \subset \overline{\mathbb{H}_\mathbb{C}^2} \backslash \{\infty\} = \mathbb{C} \times \mathbb{R} \times [0, +\infty)$ are correspondingly the cones of the $\epsilon$-neighborhoods of the points $(1, 0, 0), (-1, 0, 0) \in \mathbb{C} \times \mathbb{R} \times [0, +\infty)$ with respect to the Cygan metric $\rho_c$ in $\overline{\mathbb{H}_\mathbb{C}^2} \backslash \{\infty\}$, see (2.1).

Clearly, we may consider the length $\ell$ of the geodesic $\alpha$ so small that closures of all equidistant halfspaces in $\overline{\mathbb{H}_\mathbb{C}^2} \backslash \{\infty\}$ bounded by those bisectors (and whose interiors are disjoint from the Dirichlet polyhedron $D_z(G)$) do not intersect the co-vertical bisector whose infinity is $i\mathbb{R} \times \mathbb{R} \subset \mathbb{C} \times \mathbb{R}$. It follows from the fact [G4, Thm VII.4.0.3] that equidistant half-spaces $\mathfrak{S}_1$ and $\mathfrak{S}_2$ in $\mathbb{H}_\mathbb{C}^2$ are disjoint if and only if the half-planes $\mathfrak{S}_1 \cap \mathbb{H}_\mathbb{R}^2$ and $\mathfrak{S}_2 \cap \mathbb{H}_\mathbb{R}^2$ are disjoint, see Figures 6-9.

Now we are ready to define a quasiconformal bending deformation of the group $G$ along the geodesic $A$, which defines a bending deformation of the complex surface $M = \mathbb{H}_\mathbb{C}^2/G$ along the given closed geodesic $\alpha \subset S \subset M$.

We specify numbers $\eta$ and $\zeta$ such that $0 < \zeta < \pi/2$, $0 \leq \eta < \pi - 2\zeta$ and the intersection $C_+ \cap (\mathbb{C} \times \{0\})$ is contained in the angle $\{z \in \mathbb{C} : |\arg z| \leq \zeta\}$. Then we define a bending homeomorphism $\phi = \phi_{\eta,\zeta} : \mathbb{C} \to \mathbb{C}$ which bends the real axis $\mathbb{R} \subset \mathbb{C}$ at the origin by the angle $\eta$, see Fig. 10:

$$\phi_{\eta,\zeta}(z) = \begin{cases} z & \text{if } |\arg z| \geq \pi - \zeta \\ z \cdot \exp(i\eta) & \text{if } |\arg z| \leq \zeta \\ z \cdot \exp(i\eta(1 - (\arg z - \zeta)/(\pi - 2\zeta))) & \text{if } \zeta < \arg z < \pi - \zeta \\ z \cdot \exp(i\eta(1 + (\arg z + \zeta)/(\pi - 2\zeta))) & \text{if } \zeta - \pi < \arg z < -\zeta \,. \end{cases} \quad (5.1)$$

FIGURE 10

For negative $\eta$, $2\zeta - \pi < \eta < 0$, we set $\phi_{\eta,\zeta}(z) = \overline{\phi_{-\eta,\zeta}(\overline{z})}$. Clearly, $\phi_{\eta,\zeta}$ is quasiconformal with respect to the Cygan norm (2.1) and is an isometry in the $\zeta$-cone neighborhood of the real axis $\mathbb{R}$ because its linear distortion is given by

$$K(\phi_{\eta,\zeta}, z) = \begin{cases} 1 & \text{if } |\arg z| \geq \pi - \zeta \\ 1 & \text{if } |\arg z| \leq \zeta \\ (\pi - 2\zeta)/(\pi - 2\zeta - \eta) & \text{if } \zeta < \arg z < \pi - \zeta \\ (\pi - 2\zeta + \eta)/(\pi - 2\zeta) & \text{if } \zeta - \pi < \arg z < -\zeta \,. \end{cases} \quad (5.2)$$

Foliating the punctured Heisenberg group $\mathcal{H}\backslash\{0\}$ by Heisenberg spheres $S(0,r)$ of radii $r > 0$, we can extend the bending homeomorphism $\phi_{\eta,\zeta}$ to an elementary bending homeomorphism $\varphi = \varphi_{\eta,\zeta} : \mathcal{H} \to \mathcal{H}$, $\varphi(0) = 0$, $\varphi(\infty) = \infty$, of the whole sphere $S^3 = \overline{\mathcal{H}}$ at infinity.

Namely, for the "dihedral angles" $W_+, W_- \subset \mathcal{H}$ with the common vertical axis $\{0\} \times \mathbb{R}$ and which are foliated by arcs of real circles connecting points $(0,v)$ and $(0,-v)$ on the vertical axis and intersecting the the $\zeta$-cone neighborhoods of infinite rays $\mathbb{R}_+, \mathbb{R}_- \subset \mathbb{C}$, correspondingly, the restrictions $\varphi|_{W_-}$ and $\varphi|_{W_+}$ of the bending homeomorphism $\varphi = \varphi_{\eta,\zeta}$ are correspondingly the identity and the unitary rotation $U_\eta \in PU(2,1)$ by angle $\eta$ about the vertical axis $\{0\} \times \mathbb{R} \subset \mathcal{H}$, see also [A10, (4.4)]. Then it follows from (5.2) that $\varphi_{\eta,\zeta}$ is a $G_0$-equivariant quasiconformal homeomorphism in $\overline{\mathcal{H}}$.

We can naturally extend the foliation of the punctured Heisenberg group $\mathcal{H}\backslash\{0\}$ by Heisenberg spheres $S(0,r)$ to a foliation of the hyperbolic space $\mathbb{H}^2_\mathbb{C}$ by bisectors $\mathfrak{S}_r$ having those $S(0,r)$ as the spheres at infinity. It is well known (see [Mo2]) that each bisector $\mathfrak{S}_r$ contains a geodesic $\gamma_r$ which connects points $(0,-r^2)$ and $(0,r^2)$ of the Heisenberg group $\mathcal{H}$ at infinity, and furthermore $\mathfrak{S}_r$ fibers over $\gamma_r$ by complex geodesics $Y$ whose circles at infinity are complex circles foliating the sphere $S(0,r)$.

Using those foliations of the hyperbolic space $\mathbb{H}^2_\mathbb{C}$ and bisectors $\mathfrak{S}_r$, we extend the elementary bending homeomorphism $\varphi_{\eta,\zeta} : \overline{\mathcal{H}} \to \overline{\mathcal{H}}$ at infinity to an elementary bending homeomorphism $\Phi_{\eta,\zeta} : \overline{\mathbb{H}^2_\mathbb{C}} \to \overline{\mathbb{H}^2_\mathbb{C}}$. Namely, the map $\Phi_{\eta,\zeta}$ preserves each of bisectors $\mathfrak{S}_r$, each complex geodesic fiber $Y$ in such bisectors, and fixes the intersection points $y$ of those complex geodesic fibers and the complex geodesic connecting the origin and $\infty$ of the Heisenberg group $\mathcal{H}$ at infinity. We complete our extension $\Phi_{\eta,\zeta}$ by defining its restriction to a given (invariant) complex geodesic fiber $Y$ with the fixed point $y \in Y$. This map is obtained by radiating the circle homeomorphism $\varphi_{\eta,\zeta}|_{\partial Y}$ to the whole (Poincaré) hyperbolic 2-plane $Y$ along geodesic rays $[y,\infty) \subset Y$, so that it preserves circles in $Y$ centered at $y$ and bends (at $y$, by the angle $\eta$) the geodesic in $Y$ connecting the central points of the corresponding arcs of the complex circle $\partial Y$, see Fig.10.

Due to the construction, the elementary bending (quasiconformal) homeomorphism $\Phi_{\eta,\zeta}$ commutes with elements of the cyclic loxodromic group $G_0 \subset G$. Another most important property of the homeomorphism $\Phi_{\eta,\zeta}$ is the following.

Let $D_z(G)$ be the Dirichlet fundamental polyhedron of the group $G$ centered at a given point $z$ on the axis $A$ of the cyclic loxodromic group $G_0 \subset G$, and $\mathfrak{S}^+ \subset \mathbb{H}^2_\mathbb{C}$ be a "half-space" disjoint from $D_z(G)$ and bounded by a bisector $\mathfrak{S} \subset \mathbb{H}^2_\mathbb{C}$ which is different from bisectors $\mathfrak{S}_r, r > 0$, and contains a side $\mathfrak{s}$ of the polyhedron $D_z(G)$. Then there is an open neighborhood $U(\overline{\mathfrak{S}^+}) \subset \overline{\mathbb{H}^2_\mathbb{C}}$ such that the restriction of the elementary bending homeomorphism $\Phi_{\eta,\zeta}$ to it either is the identity or coincides with the unitary rotation $U_\eta \subset PU(2,1)$ by the angle $\eta$ about the "vertical" complex geodesic (containing the vertical axis $\{0\} \times \mathbb{R} \subset \mathcal{H}$ at infinity).

The above properties of quasiconformal homeomorphism $\Phi = \Phi_{\eta,\zeta}$ show that the image $D_\eta = \Phi_{\eta,\zeta}(D_z(G))$ is a polyhedron in $\mathbb{H}^2_\mathbb{C}$ bounded by bisectors. Furthermore, there is a natural identification of its sides induced by $\Phi_{\eta,\zeta}$. Namely, the pairs of sides preserved by $\Phi$ are identified by the original generators of the group $G_1 \subset G$. For other sides $\mathfrak{s}_\eta$ of $D_\eta$, which are images of corresponding sides $\mathfrak{s} \subset D_z(G)$ under the unitary rotation $U_\eta$, we define side pairings by using the group $G$ decomposition (see Fig. 6-9).

Actually, if $G = G_1 *_{G_0} G_2$, we change the original side pairings $g \in G_2$ of $D_z(G)$-sides to the hyperbolic isometries $U_\eta g U_\eta^{-1} \in PU(2,1)$. In the case of HNN-extension, $G = G_1 *_{G_0} = \langle G_1, g_2 \rangle$, we change the original side pairing $g_2 \in G$ of $D_z(G)$-sides to the hyperbolic isometry $U_\eta g_2 \in PU(2,1)$. In other words, we define deformed groups $G_\eta \subset PU(2,1)$ correspondingly as

$$G_\eta = G_1 *_{G_0} U_\eta G_2 U_\eta^{-1} \quad \text{or} \quad G_\eta = \langle G_1, U_\eta g_2 \rangle = G_1 *_{G_0} . \tag{5.3}$$

This shows that the family of representations $G \to G_\eta \subset PU(2,1)$ does not depend on angles $\zeta$ and holomorphically depends on the angle parameter $\eta$. Let us also observe that, for small enough angles $\eta$, the behavior of neighboring polyhedra $g'(D_\eta)$, $g' \in G_\eta$ is the same as of those $g(D_z(G))$, $g \in G$, around the Dirichlet fundamental polyhedron $D_z(G)$. This is because the new polyhedron $D_\eta \subset \mathbb{H}_{\mathbb{C}}^2$ has isometrically the same (tesselations of) neighborhoods of its side-intersections as $D_z(G)$ had. This implies that the polyhedra $g'(D_\eta)$, $g' \in G_\eta$, form a tesselation of $\mathbb{H}_{\mathbb{C}}^2$ (with non-overlapping interiors). Hence the deformed group $G_\eta \subset PU(2,1)$ is a discrete group, and $D_\eta$ is its fundamental polyhedron bounded by bisectors.

Using $G$-compatibility of the restriction of the elementary bending homeomorphism $\Phi = \Phi_{\eta,\zeta}$ to the closure $\overline{D_z(G)} \subset \overline{\mathbb{H}_{\mathbb{C}}^2}$, we equivariantly extend it from the polyhedron $\overline{D_z(G)}$ to the whole space $\mathbb{H}_{\mathbb{C}}^2 \cup \Omega(G)$ accordingly to the $G$-action. In fact, in terms of the natural isomorphism $\chi : G \to G_\eta$ which is identical on the subgroup $G_1 \subset G$, we can write the obtained $G$-equivariant homeomorphism $F = F_\eta : \overline{\mathbb{H}_{\mathbb{C}}^2} \setminus \Lambda(G) \to \overline{\mathbb{H}_{\mathbb{C}}^2} \setminus \Lambda(G_\eta)$ in the following form:

$$\begin{aligned} F_\eta(x) &= \Phi_\eta(x) \quad \text{for} \quad x \in \overline{D_z(G)}, \\ F_\eta \circ g(x) &= g_\eta \circ F_\eta(x) \quad \text{for} \quad x \in \overline{\mathbb{H}_{\mathbb{C}}^2} \setminus \Lambda(G),\ g \in G,\ g_\eta = \chi(g) \in G_\eta . \end{aligned} \tag{5.4}$$

Due to quasiconformality of $\Phi_\eta$, the extended $G$-equivariant homeomorphism $F_\eta$ is quasiconformal. Furthermore, its extension by continuity to the limit (real) circle $\Lambda(G)$ coincides with the canonical equivariant homeomorphism $f_\chi : \Lambda(G) \to \Lambda(G_\eta)$ given by the isomorphism Theorem 3.2. Hence we have a $G$-equivariant quasiconformal self-homeomorphism of the whole space $\overline{\mathbb{H}_{\mathbb{C}}^2}$, which we denote as before by $F_\eta$.

The family of $G$-equivariant quasiconformal homeomorphisms $F_\eta$ induces representations $F_\eta^* : G \to G_\eta = F_\eta G_2 F_\eta^{-1}$, $\eta \in (-\eta_0, \eta_0)$. In other words, we have a curve $\mathcal{B} : (-\eta_0, \eta_0) \to \mathcal{R}_0(G)$ in the variety $\mathcal{R}_0(G)$ of faithful discrete representations of $G$ into $PU(2,1)$, which covers a nontrivial curve in the Teichmüller space $\mathcal{T}(G)$ represented by conjugacy classes $[\mathcal{B}(\eta)] = [F_\eta^*]$. We call the constructed deformation $\mathcal{B}$ the bending deformation of a given lattice $G \subset PO(2,1) \subset PU(2,1)$ along a bending geodesic $A \subset \mathbb{H}_{\mathbb{C}}^2$ with loxodromic stabilizer $G_0 \subset G$. In terms of manifolds, $\mathcal{B}$ is the bending deformation of a given complex surface $M = \mathbb{H}_{\mathbb{C}}^2/G$ homotopy equivalent to its totally real geodesic surface $S_g \subset M$, along a given simple geodesic $\alpha$.

□

*Remark 5.4.* It follows from the above construction of the bending homeomorphism $F_{\eta,\zeta}$, that the deformed complex hyperbolic surface $M_\eta = \mathbb{H}_{\mathbb{C}}^2/G_\eta$ fibers over the pleated hyperbolic surface $S_\eta = F_\eta(\mathbb{H}_{\mathbb{R}}^2)/G_\eta$ (with the closed geodesic $\alpha$ as the singular locus). The fibers of this fibration are "singular real planes" obtained from totally real geodesic 2-planes by bending them by angle $\eta$ along complete real geodesics. These (singular) real geodesics are the intersections of the complex geodesic connecting the axis $A$ of the cyclic group $G_0 \subset G$ and the totally real geodesic planes that represent fibers of the original fibration in $M = \mathbb{H}^2_\mathbb{C}/G$.

*Proof of Corollary 5.2.* Since, due to (5.3), bendings along disjoint closed geodesics are independent, we need to show that our bending deformation is not trivial, and $[\mathcal{B}(\eta)] \neq [\mathcal{B}(\eta')]$ for any $\eta \neq \eta'$.

The non-triviality of our deformation follows directly from (5.3), compare [A9]. Namely, the restrictions $\rho_\eta|_{G_1}$ of bending representations to a non-elementary subgroup $G_1 \subset G$ (in general, to a "real" subgroup $G_r \subset G$ corresponding to a totally real geodesic piece in the homotopy equivalent surface $S \simeq M$) are identical. So if the deformation $\mathcal{B}$ were trivial then it would be conjugation of the group $G$ by projective transformations that commute with the non-trivial real subgroup $G_r \subset G$ and pointwise fix the totally real geodesic plane $\mathbb{H}^2_\mathbb{R}$. This contradicts to the fact that the limit set of any deformed group $G_\eta$, $\eta \neq 0$, does not belong to the real circle containing the limit Cantor set $\Lambda(G_r)$.

The injectivity of the map $\mathcal{B}$ can be obtained by using Élie Cartan [Car] angular invariant $\mathbb{A}(x)$, $-\pi/2 \leq \mathbb{A}(x) \leq \pi/2$, for a triple $x = (x^0, x^1, x^2)$ of points in $\partial \mathbb{H}^2_\mathbb{C}$. It is known (see [G4]) that, for two triples $x$ and $y$, $\mathbb{A}(x) = \mathbb{A}(y)$ if and only if there exists $g \in PU(2,1)$ such that $y = g(x)$; furthermore, such a $g$ is unique provided that $\mathbb{A}(x)$ is neither zero nor $\pm \pi/2$. Here $\mathbb{A}(x) = 0$ if and only if $x^0, x^1$ and $x^2$ lie on an $\mathbb{R}$-circle, and $\mathbb{A}(x) = \pm \pi/2$ if and only if $x^0, x^1$ and $x^2$ lie on a chain ($\mathbb{C}$-circle).

Namely, let $g_2 \in G \backslash G_1$ be a generator of the group $G$ in (4.5) whose fixed point $x^2 \in \Lambda(G)$ lies in $\mathbb{R}_+ \times \{0\} \subset \mathcal{H}$, and $x^2_\eta \in \Lambda(G_\eta)$ the corresponding fixed point of the element $\chi_\eta(g_2) \in G_\eta$ under the free-product isomorphism $\chi_\eta : G \to G_\eta$. Due to our construction, one can see that the orbit $\gamma(x^2_\eta)$, $\gamma \in G_0$, under the loxodromic (dilation) subgroup $G_0 \subset G \cap G_\eta$ approximates the origin along a ray $(0, \infty)$ which has a non-zero angle $\eta$ with the ray $\mathbb{R}_- \times \{0\} \subset \mathcal{H}$. The latter ray also contains an orbit $\gamma(x^1)$, $\gamma \in G_0$, of a limit point $x^1$ of $G_1$ which approximates the origin from the other side. Taking triples $x = (x^1, 0, x^2)$ and $x_\eta = (x^1, 0, x^2_\eta)$ of points which lie correspondingly in the limit sets $\Lambda(G)$ and $\Lambda(G_\eta)$, we have that $\mathbb{A}(x) = 0$ and $\mathbb{A}(x_\eta) \neq 0, \pm \pi/2$. Due to Theorem 3.2, both limit sets are topological circles which however cannot be equivalent under a hyperbolic isometry because of different Cartan invariants (and hence, again, our deformation is not trivial).

Similarly, for two different values $\eta$ and $\eta'$, we have triples $x_\eta$ and $x_{\eta'}$ with different (non-trivial) Cartan angular invariants $\mathbb{A}(x_\eta) \neq \mathbb{A}(x_{\eta'})$. Hence $\Lambda(G_\eta)$ and $\Lambda(G_{\eta'})$ are not $PU(2,1)$-equivalent.

□

One can apply the above proof to a general situation of bending deformations of a complex hyperbolic surface $M = \mathbb{H}^2_\mathbb{C}/G$ whose holonomy group $G \subset PU(2,1)$ has a non-elementary subgroup $G_r$ preserving a totally real geodesic plane $\mathbb{H}^2_\mathbb{R}$. In other words, such a complex surfaces $M$ has an embedded totally real geodesic surface with geodesic boundary. So we immediately have:

**Corollary 5.5.** *Let $M = \mathbb{H}^2_\mathbb{C}/G$ be a complex hyperbolic surface with embedded totally real geodesic surface $S_r \subset M$ with geodesic boundary, and $\mathcal{B} : (-\eta, \eta) \to \mathcal{D}(M)$ be the bending deformation of $M$ along a simple closed geodesic $\alpha \subset S_r$. Then the map $\pi \circ \mathcal{B} : (-\eta, \eta) \to \mathcal{T}(M) = \mathcal{D}(M)/PU(2,1)$ is a smooth embedding* *provided that the limit set $\Lambda(G)$ of the holonomy group $G$ does not belong to the $G$-orbit of the real circle $S_{\mathbb{R}}^1$ and the chain $S_{\mathbb{C}}^1$, where the latter is the infinity of the complex geodesic containing a lift $\tilde{\alpha} \subset \mathbb{H}_{\mathbb{C}}^2$ of the closed geodesic $\alpha$, and the former one contains the limit set of the holonomy group $G_r \subset G$ of the geodesic surface $S_r$.* □

As an application of the constructed bending deformations, we answer a well known question about cusp groups on the boundary of the Teichmüller space $\mathcal{T}(M)$ of a Stein complex hyperbolic surface $M$ fibering over a compact Riemann surface of genus $p > 1$. It is a direct corollary of the following result, see [AG]:

**Theorem 5.6.** *Let $G \subset PO(2,1) \subset PU(2,1)$ be a uniform lattice isomorphic to the fundamental group of a closed surface $S_p$ of genus $p \geq 2$. Then, for any simple closed geodesic $\alpha \subset S_p = H_{\mathbb{R}}^2/G$, there is a continuous deformation $\rho_t = f_t^*$ induced by $G$-equivariant quasiconformal homeomorphisms $f_t : \overline{\mathbb{H}_{\mathbb{C}}^2} \to \overline{\mathbb{H}_{\mathbb{C}}^2}$ whose limit representation $\rho_\infty$ corresponds to a boundary cusp point of the Teichmüller space $\mathcal{T}(G)$, that is the boundary group $\rho_\infty(G)$ has an accidental parabolic element $\rho_\infty(g_\alpha)$ where $g_\alpha \in G$ represents the geodesic $\alpha \subset S_p$.*

We note that, due to our construction of such continuous quasiconformal deformations in [AG], they are independent if the corresponding geodesics $\alpha_i \subset S_p$ are disjoint. It implies the existence of a boundary group in $\partial \mathcal{T}(G)$ with "maximal" number of non-conjugate accidental parabolic subgroups:

**Corollary 5.7.** *Let $G \subset PO(2,1) \subset PU(2,1)$ be a uniform lattice isomorphic to the fundamental group of a closed surface $S_p$ of genus $p \geq 2$. Then there is a continuous deformation $R : \mathbb{R}^{3p-3} \to \mathcal{T}(G)$ whose boundary group $G_\infty = R(\infty)(G)$ has $(3p-3)$ non-conjugate accidental parabolic subgroups.*

Department of Mathematics, University of Oklahoma, Norman, OK 73019
*E-mail address*: `apanasov@ou.edu`

Mathematical Sciences Research Institute, Berkeley, CA 94720-5070

Grad. School of Math. Sci., Univ. of Tokyo, 3-8-1 Komaba, Tokyo 153